\documentclass[11pt]{amsart}
\usepackage{amsfonts,amsmath,amssymb}%,mathrsfs}
\usepackage[dvips]{graphicx} 
\usepackage{psfrag}

\newdimen\margin   % needed for macros \textdisplay & \ltextdisplay
\def\COMMENT#1{}
\def\TASK#1{}
\def\noproof{{\unskip\nobreak\hfill\penalty50\hskip2em\hbox{}\nobreak\hfill%
       $\square$\parfillskip=0pt\finalhyphendemerits=0\par}\goodbreak}
\def\endproof{\noproof\bigskip}
\def\proof{\removelastskip\penalty55\medskip\noindent{\bf Proof. }}
\def\textno#1&#2\par{%
   \margin=\hsize
   \advance\margin by -4\parindent
          \setbox1=\hbox{\sl#1}%
   \ifdim\wd1 < \margin
      $$\box1\eqno#2$$%
   \else
      \bigbreak
      \hbox to \hsize{\indent$\vcenter{\advance\hsize by -3\parindent
      \sl\noindent#1}\hfil#2$}%
      \bigbreak
   \fi}
\def\enddiscard{}
\long\def\discard#1\enddiscard{}

\newcommand{\eps}{\varepsilon}
\newcommand{\Nat}{\mathbb{N}}
\newcommand{\Hy}{\mathcal{H}}

\newcommand{\G}{{\mathcal G}}

\newcommand{\cn}{{\mathcal N}_h}
\newcommand{\F}{{\mathcal F}}
\newcommand{\N}{{\mathcal N}}
\newcommand{\B}{{\mathcal B}}

\newcommand{\C}{\mathcal{C}}
\newcommand{\K}{\mathcal{K}}

\newcommand{\cS}{\mathcal{S}}

\newcommand{\R}{\mathcal{R}}
\newcommand{\D}{\mathcal{D}}

\newcommand{\htohp}{|H \rightarrow \Hy'|}
\newcommand{\exphtohp}{\overline{|\Hy \rightarrow \Hy'|}}
\newcommand{\htohs}{|H \rightarrow \Hy^*|}
\newcommand{\exphtohs}{\overline{|\Hy \rightarrow \Hy^*|}}
\newcommand{\nhtob}{|N_h \rightarrow \B|}
\newcommand{\expnhtob}{\overline{|\N_h \rightarrow \B|}}

\newcommand{\nhnhstofmin}{|N_h \cup N_h^*\stackrel{k-1}{\rightarrow} \F|}
\newcommand{\expnhnhstofmin}{\overline{|\N_h \cup \N_h^*\stackrel{k-1}{\rightarrow} \F|}}
\newcommand{\nhnhstof}{|N_h \cup N_h^*\rightarrow \F|}
\newcommand{\expnhnhstof}{\overline{|\N_h \cup \N_h^*\rightarrow \F|}}
\newcommand{\nhtohh}{|N_h \rightarrow \Hy_h|}

\newcommand{\usef}{\mathtt{Usef}}
\newcommand{\typ}{\mathtt{typ}}
\newcommand{\atyp}{\mathtt{atyp}}
\newcommand{\nhstohhs}{|N_h^* \rightarrow \Hy_h^*|}

\newcommand{\hhtoh}{|H_h \rightarrow \Hy|}
\newcommand{\exphhtoh}{\overline{|\Hy_h \rightarrow \Hy|}}
\newcommand{\hptohpp}{|H'\rightarrow \Hy''|}
\newcommand{\exphptohpp}{\overline{|\Hy'\rightarrow \Hy''|}}
\newcommand{\hptohppmin}{|H'\stackrel{k-1}{\rightarrow}\Hy''|}
\newcommand{\exphptohppmin}{\overline{|\Hy'\stackrel{k-1}{\rightarrow}\Hy''|}}
\newcommand{\Part}{{\mathcal P}}
\newcommand{\cross}{\mathrm {Cross}}
\newcommand{\av}{\mbox{\boldmath $a$\unboldmath}}

\newtheorem{firsttheorem}{Proposition}

\newtheorem{fact}[firsttheorem]{Fact}
\newtheorem{theorem}[firsttheorem]{Theorem}
\newtheorem{lemma}[firsttheorem]{Lemma}

\newtheorem{proposition}[firsttheorem]{Proposition}

\begin{document}
\title{Embeddings and Ramsey numbers of sparse $k$-uniform hypergraphs}
\author{Oliver Cooley, Nikolaos Fountoulakis, Daniela K\"uhn and Deryk
Osthus}
\thanks {N.~Fountoulakis and D.~K\"uhn were supported by the EPSRC, grant no.~EP/D50564X/1}
\maketitle
\begin{abstract} \noindent
Chv\'atal, R\"odl, Szemer\'edi and Trotter~\cite{CRST} proved that
the Ramsey numbers of graphs of bounded maximum degree are linear
in their order. In \cite{CNKO, NORS} the same result was proved
for $3$-uniform hypergraphs. Here we extend this result
to $k$-uniform hypergraphs for any integer $k \geq 3$. As in the
$3$-uniform case, the main new tool which we prove and use is an
embedding lemma for $k$-uniform hypergraphs of bounded maximum
degree into suitable $k$-uniform `quasi-random' hypergraphs.
\end{abstract}
\noindent
{\footnotesize
keywords: hypergraphs; regularity lemma; Ramsey numbers; embedding problems}

\vspace{0.4cm}

\section{Introduction}\label{intro}

%$\mathscr{P}$

The \emph{Ramsey number $R(\Hy)$} of a $k$-uniform hypergraph
$\Hy$ is the smallest $N \in \mathbb{N}$ such that for every
$2$-colouring of the hyperedges of the complete $k$-uniform
hypergraph on $N$ vertices one can find a monochromatic copy of
$\Hy$. For general $\Hy$, the best upper bound is due to Erd\H{o}s
and Rado~\cite{ER52}. Writing $|\Hy |$ for the number of vertices of $\Hy$,
it implies that for any $k\geq 2$
$$R(\Hy) \leq 2^{2^{\cdot^{\cdot^{\cdot^{2^{c_k|\Hy|}}}}}},$$
where the number of 2's is~$k-1$. In the other direction, Erd\H{o}s and Hajnal (see~\cite{GSR})
showed that if $k\geq 3$ and~$\Hy$ is a complete $k$-uniform hypergraph then $R(\Hy)$ is bounded below
by a tower in which the number of 2's is~$k-2$ and the top exponent is~$c_k'|\Hy|^2$.

For the case of graphs (i.e.,~when
$k=2$) it is known that there are many families of graphs $H$ for
which the Ramsey number is much smaller than exponential. In
particular, Burr and Erd\H{o}s~\cite{BE75} asked for which
graphs~$H$ the Ramsey number $R(H)$ is linear in the order $|H|$ of~$H$ and
conjectured this to be true for graphs of bounded maximum degree.
This was proved by Chv\'atal, R\"odl, Szemer\'edi and
Trotter~\cite{CRST}. Here we show that their result extends to $k$-uniform
hypergraphs~$\Hy$ of bounded maximum degree, where the \emph{degree of a vertex~$x$}
in~$\Hy$ is defined to be the number of hyperedges which contain~$x$.

\begin{theorem} \label{Ramseythm}
For all $\Delta,k\in\mathbb{N}$ there exists a
constant~$C=C(\Delta,k)$ such that all $k$-uniform hypergraphs~$\Hy$
of maximum degree at most~$\Delta$ satisfy $R(\Hy) \le C |\Hy|$.
\end{theorem}
The overall strategy of our proof of Theorem~\ref{Ramseythm} is
related to that of Chv\'atal et al.~\cite{CRST}, which is based on
the regularity lemma for graphs. We apply a version (due to R\"odl and Schacht~\cite{roedlschacht})
of the regularity lemma for $k$-uniform hypergraphs.
Roughly speaking, it guarantees a partition of an arbitrary dense $k$-uniform
hypergraph into `quasi-random' subhypergraphs. Our main contribution
is an embedding result (Theorem~\ref{embhgph}) which guarantees the existence of a copy of
a hypergraph~$\Hy$ of bounded maximum degree inside a suitable
`quasi-random' hypergraph~$\G$ even if the order of~$\Hy$ is linear
in that of~$\G$. In fact, we prove a stronger embedding result of independent
interest (Theorem~\ref{emblemma}). It even counts the number of copies of such $\Hy$ in $\G$ and
thus generalizes the well-known hypergraph counting lemma (which only allows
for bounded size $\Hy$).

After the submission of this paper, Keevash~\cite{Keev} extended
Theorem~\ref{embhgph} to a hypergraph blow-up lemma for embeddings
of spanning subhypergraphs~$\Hy$. The case of $3$-uniform
hypergraphs in Theorem~\ref{Ramseythm} was proved recently
in~\cite{CNKO} and independently by~Nagle, Olsen, R\"odl and
Schacht~\cite{NORS}. Also, Kostochka and R\"odl~\cite{RK_OW}
earlier proved an approximate version of Theorem~\ref{Ramseythm}:
for all $\varepsilon,\Delta,k>0$ there is a constant $C$ such that
$R(\Hy) \le C|\Hy|^{1+\varepsilon}$ if~$\Hy$ has maximum degree at
most~$\Delta$. After this manuscript was submitted, Conlon, Fox
and Sudakov~\cite{CFS} obtained a proof of Theorem~\ref{Ramseythm}
which does not rely on hypergraph regularity and gives a better
bound on $C$. Also, Ishigami~\cite{Ish} independently announced a
proof of Theorem~\ref{Ramseythm} using a similar approach to ours.
Apart from these, the only previous results on the Ramsey numbers
of sparse hypergraphs are on hypergraph cycles (see
e.g.~\cite{bergecycles,haxell1, haxell2}.

It would  be desirable to extend Theorem~\ref{Ramseythm} to a larger
class of hypergraphs. For instance the graph analogue of
Theorem~\ref{Ramseythm} is known for so-called $p$-arrangeable
graphs~\cite{chenschelp}, which include the class of all planar
graphs. However, R\"odl and Kostochka~\cite{RK_OW} showed that a
natural hypergraph analogue of the famous Burr-Erd\H{o}s
conjecture on Ramsey numbers of $d$-degenerate graphs fails for
$k$-uniform hypergraphs if~$k\ge 3$. (A graph is $d$-degenerate if
the maximum average degree over all its subgraphs is at most $d$.
If a graph is $p$-arrangeable, then it is also $d$-degenerate for
some $d$.) But it may still be possible to generalize the Burr-Erd\H{o}s
conjecture to hypergraphs in a different way.

This paper is organized as follows.
In Section~\ref{sec:emblemma} we give an overview of the proof of
Theorem~\ref{Ramseythm} and we state the embedding theorem (Theorem~\ref{embhgph})
mentioned above. Our proof of Theorem~\ref{embhgph} relies on a more general version
(Lemma~\ref{countinglemma}) of the well-known counting
lemma for hypergraphs as well as an `extension lemma' (Lemma~\ref{extensions_count}),
whose proofs are postponed until Sections~\ref{sec:CL} and~\ref{sec:ext}.
We introduce these lemmas, along with further tools, in Section~\ref{sec:tools}.
We then prove a strengthened version (Theorem~\ref{emblemma}) of Theorem~\ref{embhgph} in
Section~\ref{sec:embed}. The regularity lemma for $k$-uniform
hypergraphs is introduced in Section~\ref{regularity}. In
Section~\ref{sec:proofofThm1} we deduce Theorem~\ref{Ramseythm} from the
regularity lemma and Theorem~\ref{embhgph}.
In Section~\ref{sec:CL} we derive our version of the counting lemma (Lemma~\ref{countinglemma})
from that in~\cite{roedlschacht2}. Finally, in Section~\ref{sec:ext} we use it to deduce the
extension lemma (Lemma~\ref{extensions_count}).

\section{Overview of the proof of Theorem~\ref{Ramseythm} and statement
of the embedding theorem}\label{sec:emblemma}

\subsection{Overview of the proof of Theorem~\ref{Ramseythm}}
The proof in \cite{CRST} that graphs of bounded degree have linear
Ramsey numbers proceeds roughly as follows: Let $H$ be a graph of
maximum degree $\Delta$. Take a complete graph $K_n$, where $n$ is
a sufficiently large integer. Colour the edges of $K_n$ with red and blue, and
apply the graph regularity lemma to the denser of the two
monochromatic graphs, $G_{red}$ say, to obtain a partition of the
vertex set into a bounded number of clusters. Since almost all pairs
of clusters are regular or `quasi-random', by Tur\'an's theorem
there will be a set of $r$ clusters, where $r:=R(K_{\Delta+1})$,
in which each pair of clusters is regular. A pair of clusters will be
coloured red if its density in~$G_{red}$ is at least~$1/2$, and blue
otherwise. By the definition of~$r$, there must be a set
of~$\Delta+1$ clusters such that all the pairs have the same
colour. If this colour is red, then one can apply the so-called
embedding or key lemma for graphs to find a (red) copy of~$H$ in
the subgraph of~$G_{red}$ spanned by these $\Delta+1$ clusters. This
is possible since $\chi(H) \leq \Delta+1$. If all the pairs of
clusters are coloured blue we apply the embedding theorem in the blue
subgraph $G_{blue}$ of~$K_n$ to find a blue copy of~$H$. It turns out that in this
proof we only needed $n \geq C|H|$, where $C$ is a
constant dependent only on~$\Delta$. Thus $R(H) \leq C|H|$.

We will generalize this approach to $k$-uniform hypergraphs. As mentioned in
Section~\ref{intro}, the main obstacle is the proof of
an embedding theorem for $k$-uniform hypergraphs (Theorem~\ref{embhgph} below), which
allows us to embed a $k$-uniform hypergraph~$\Hy$ within a suitable `quasi-random'
$k$-uniform hypergraph~$\G$, where the order of~$\Hy$ might be linear
in the order of~$\G$. Our proof uses ideas from~\cite{CNKO}.

\subsection{Notation and statement of the embedding theorem}
Before we can state the embedding theorem, we first have to say
what we mean by a regular or `quasi-random' hypergraph. In the
setup below, this will involve the relationship between certain
$i$-uniform hypergraphs and $(i-1)$-uniform hypergraphs on the
same vertex set. Given a hypergraph $\G$, we write $E(\G)$ for the
set of its hyperedges and define $e(\G):=|E(\G)|$. We
write~$K^{(j)}_i$ for the complete $j$-uniform hypergraph on~$i$
vertices. Given a $j$-uniform hypergraph~$\G$ and
$j \leq i$, we write $\K_i(\G)$ for the set of $i$-sets
of vertices of~$\G$ which form a copy of~$K_i^{(j)}$ in~$\G$.
Given an $i$-partite $i$-uniform hypergraph~$\G_i$, and an
$i$-partite $(i-1)$-uniform hypergraph~$\G_{i-1}$ on the same
vertex set, we
define the \emph{density of $\G_i$ with respect to $\G_{i-1}$} to be%
$$
d(\G_i|\G_{i-1}):= \frac{|\K_i(\G_{i-1})\cap E(\G_i)|}{|\K_i(\G_{i-1})|}
$$
if $|\K_i(\G_{i-1})|>0$, and $d(\G_i|\G_{i-1}):=0$ otherwise.
More generally, if $\mathbf{Q}:=(Q(1),Q(2),\ldots,Q(r))$ is a
collection of~$r$ subhypergraphs of~$\G_{i-1}$, we define
$\K_i(\mathbf{Q}):=\bigcup_{j=1}^r \K_i(Q(j))$ and
$$
d(\G_i |\mathbf{Q}):= \frac{|\K_i(\mathbf{Q})\cap E(\G_i)|}{|\K_i(\mathbf{Q})|}
$$
if $|\K_i(\mathbf{Q})|>0$, and $d(\G_i |\mathbf{Q}):=0$ otherwise.
We sometimes write $|K^{(j)}_i|_\mathbf{Q}$ instead of
$|\K_i(\mathbf{Q})|$.

We say that $\G_i$ is \emph{$(d_i,\delta,r)$-regular with respect
to~$\G_{i-1}$} if every $r$-tuple $\mathbf{Q}$ with
$|\K_i(\mathbf{Q})|>\delta|\K_i( \G_{i-1})|$ satisfies
$$
d(\G_i|\mathbf{Q})= d_i \pm \delta.
$$
Given $\ell\ge i\ge 3$, an $\ell$-partite $i$-uniform hypergraph~$\G_i$ and an $\ell$-partite
$(i-1)$-uniform hypergraph~$\G_{i-1}$ on the same vertex set, we say that
$\G_i$ is $(d_i,\delta,r)$-regular with respect to~$\G_{i-1}$ if
for every $i$-tuple~$K$ of vertex classes, either~$\G_i[K]$ is
$(d_i,\delta,r)$-regular with respect to~$\G_{i-1}[K]$ or
$d(\G_i[K] | \G_{i-1}[K])=0$ (but the latter should not hold for
all~$K$).%
     \COMMENT{The last part ensures that this def agrees with the one
for $\ell=i$.}
Instead of $(d_i,\delta,1)$-regularity
we sometimes refer to \emph{$(d_i,\delta)$-regularity}.

The density of a bipartite graph~$G$ with vertex classes~$A$
and~$B$ is defined by $d(A,B):=e(A,B)/|A||B|$ and~$G$ is $(d,\delta)$-regular if
for all sets $X\subseteq A$ and $Y\subseteq B$ with $|X|\ge \delta|A|$ and $|Y|\ge \delta|B|$
we have $d(X,Y)=d\pm \delta$.
We say that an $\ell$-partite graph~$\G_2$ is \emph{$(d_2,\delta)$-regular} if each
of the $\binom{\ell}{2}$ bipartite subgraphs forming it is either $(d_2,\delta)$-regular
or has density~$0$ (and if for at least one of them the former holds).

Suppose that we have $\ell \geq k$ vertex classes $V_1,\ldots,V_\ell$, and
that for each $i=2,\dots,k$ we are given an $\ell$-partite $i$-uniform
hypergraph~$\G_i$ with these vertex classes. Suppose also that $\Hy$ is
an $\ell$-partite $k$-uniform hypergraph with vertex classes
$X_1,\ldots,X_\ell$. We will aim to embed~$\Hy$ into~$\G_k$, and in particular
to embed~$X_j$ into~$V_j$ for each $j=1,\dots,\ell$. So we make the following definition:
We say that $(\G_k,\ldots,\G_2)$ \emph{respects the partition of~$\Hy$}
if whenever~$\Hy$ contains a hyperedge with vertices in
$X_{j_1},\ldots,X_{j_k}$, then there is a hyperedge of~$\G_k$ with
vertices in $V_{j_1},\ldots,V_{j_k}$ which also forms a copy of~$K_k^{(i)}$
in~$\G_i$ for each $i=2,\dots,k-1$.

\begin{theorem}[Embedding theorem for hypergraphs]\label{embhgph}
Let $\Delta,k,\ell,r,n_0$ be positive integers with $k\le \ell$ and let
$c,d_2,d_3,\ldots,d_k,\delta,\delta_k$ be positive constants such that
$1/d_i \in\mathbb{N}$,
$$
1/n_0\ll 1/r,\delta \ll\min\{\delta_k,d_2,\ldots,d_{k-1}\}\le\delta_k\ll
d_k,1/\Delta,1/\ell
$$
and
$$c\ll d_2,\ldots,d_k,1/\Delta,1/\ell.$$
Then the following holds for all integers $n\ge n_0$. Suppose that
$\Hy$ is an $\ell$-partite $k$-uniform hypergraph of maximum
degree at most~$\Delta$ with vertex classes $X_1,\dots,X_\ell$
such that $|X_i|\le cn$ for all $i=1,\dots,\ell$. Suppose that for
each $i=2,\dots,k$, $\G_{i}$ is an $\ell$-partite $i$-uniform
hypergraph with vertex classes $V_1,\dots,V_\ell$, which all have
size~$n$. Suppose also that~$\G_{k}$ is
$(d_{k},\delta_{k},r)$-regular with respect to $\G_{k-1}$, that
for each $i=3,\dots,k-1$, $\G_{i}$ is $(d_{i},\delta)$-regular
with respect to $\G_{i-1}$, that $\G_2$ is $(d_2,\delta)$-regular,
and that $(\G_k,\dots,\G_2)$ respects the partition of~$\Hy$.
Then~$\G_k$ contains a copy of~$\Hy$.
\end{theorem}

In the statement of Theorem~\ref{embhgph} we used the following
notation (which will be used frequently later on as well). Given
constants $a_1,a_2,a_3$, we write $a_1 \ll a_2 \ll a_3$ to mean
that we choose these constants from right to left, and there are
increasing functions $f$ and $g$ such that the lemma holds
provided that $a_2 \leq f(a_3)$ and $a_1 \leq g(a_2)$. The
functions $f$ and $g$ are determined by the calculations in the
proof of Theorem~\ref{embhgph}, but for clarity of the exposition we
will not determine them explicitly.%
     \COMMENT{The relationship between $1/r$ and $\delta$ is not made
explicit in the Reg lemma or the counting lemma of R\"odl and Schacht.
However, if $r'$ and
$\delta'$ are the constants which work in these lemmas, then we can
interpret our $1/r,\delta\ll ...$ as $1/r,\delta\le \min\{1/r',\delta'\}$.
Then our hypergraphs will be regular enough for these lemmas.}

\section{Further notation and tools}\label{sec:tools}

\subsection{Embedding theorem for complexes}
Instead of Theorem~\ref{embhgph}, we will prove a considerably stronger version
which appears as Theorem~\ref{emblemma} below. It allows the embedding of hypergraphs
which are not necessarily uniform and gives a lower bound on the number of such
embeddings. This enables us to prove the lemma by induction on~$|\Hy|$. Before
we can state Theorem~\ref{emblemma}, we need to make the following definitions.

A \emph{complex~$\Hy$} on a vertex set~$V$ is a collection of
subsets of~$V$, each of size at least~$2$, such that if $B \in
\Hy$, and if $A \subseteq B$ has size at least~$2$, then $A \in \Hy$.
(So if we add each vertex in~$V$ as a singleton into a complex, we obtain a downset.)
A \emph{$k$-complex} is a complex in which no member has
size greater than~$k$. The members of size~$i\ge 2$ are called the
\emph{$i$-edges} of~$\Hy$ and the elements of~$V$ are called the \emph{vertices}
of~$\Hy$. We write~$E_i(\Hy)$ for the set of all $i$-edges of~$\Hy$ and
set $e_i(\Hy):=|E_i(\Hy)|$. We also write
$|\Hy|:=|V|$ for the \emph{order} of~$\Hy$. Note that a $k$-uniform hypergraph can be
turned into a $k$-complex by making every hyperedge into a
complete $i$-uniform hypergraph~$K_k^{(i)}$, for each $2 \leq i\leq k$.
(In a more general $k$-complex we may have $i$-edges which do
not lie within an $(i+1)$-edge.) Given $k\le \ell$, a \emph{$(k,\ell)$-complex} is an
$\ell$-partite $k$-complex. Given a $k$-complex $\Hy$, for each
$i=2,\ldots,k$ we denote by $\Hy_i$ the \emph{underlying $i$-uniform
hypergraph of~$\Hy$}. So the vertices of~$\Hy_i$ are those of~$\Hy$ and the
hyperedges of~$\Hy_i$ are the $i$-edges of~$\Hy$.

Two vertices~$x$ and~$y$ in a $k$-complex are
\emph{neighbours} if they are joined by a $2$-edge. (Note that
if~$x$ and~$y$ lie in a common $i$-edge for some $2 \leq i \leq k$, then
they are joined by a $2$-edge.) The \emph{degree~$d(x)$} of a vertex~$x$
is the maximum (over $2\leq i \leq k$) of the number of $i$-edges containing~$x$.
Thus~$x$ has at most $d(x)$ neighbours. The \emph{maximum degree} of the
complex~$\Hy$ is the greatest degree of any vertex. Note that if~$\Hy$
is a $k$-uniform hypergraph of maximum degree~$\Delta$,  the maximum
degree of the corresponding $k$-complex is crudely at most $\Delta 2^k$.
The \emph{distance} between two vertices~$x$ and~$y$ in a $k$-complex~$\Hy$
is the length of the shortest path between~$x$ and~$y$ in the underlying
$2$-graph~$\Hy_2$ of~$\Hy$. The \emph{components} of~$\Hy$ are the subcomplexes
induced by the components of~$\Hy_2$.

We say that a $k$-complex~$\G$ is $(d_k,\ldots,d_2,\delta_k,\delta,r)$-regular
if~$\G_k$ is $(d_k,\delta_k,r)$-regular with respect to~$\G_{k-1}$, if $\G_i$
is $(d_i,\delta)$-regular with respect to $\G_{i-1}$ for each
$i=3,\dots,k-1$, and if~$\G_2$ is $(d_2,\delta)$-regular. We denote
$(d_k,\ldots,d_2)$ by $\mathbf{d}$ and refer to
$(\mathbf{d},\delta_k,\delta,r)$-regularity.

Suppose that~$\G$ is a $(k,\ell)$-complex with vertex classes $V_1,\dots,V_\ell$,
which all have size~$n$. Suppose also that~$\Hy$ is a $(k,\ell)$-complex with
vertex classes $X_1,\dots,X_\ell$ of size at most~$n$.
Similarly as for hypergraphs we say that~$\G$ \emph{respects the partition of~$\Hy$}
if whenever~$\Hy$ contains an $i$-edge with vertices in
$X_{j_1},\ldots,X_{j_i}$, then there is an $i$-edge of~$\G$ with
vertices in $V_{j_1},\ldots,V_{j_i}$.
On the other hand, we say that a labelled copy of~$\Hy$ in~$\G$ is
\emph{partition-respecting} if for each $i=1,\dots,\ell$ the vertices
corresponding to those in~$X_i$ lie within~$V_i$.
We denote by~$|\Hy|_\G$ the number of labelled, partition-respecting copies
of~$\Hy$ in~$\G$.

\begin{theorem}[Embedding theorem for complexes]\label{emblemma}
Let $\Delta,k,\ell,r,n_0$ be positive integers and let
$c,\alpha,d_2,\ldots,d_k,\delta,\delta_k$ be positive
constants such that $1/d_i \in\mathbb{N}$,
$$1/n_0\ll 1/r, \delta \ll\min\{\delta_k,d_2,\ldots,d_{k-1}\}\le\delta_k\ll
\alpha \ll d_k,1/\Delta,1/\ell$$
and
$$c\ll \alpha, d_2,\ldots,d_k.$$%
Then the following holds for all integers $n\ge n_0$. Suppose that
$\Hy$ is a $(k,\ell)$-complex of maximum degree at most~$\Delta$
with vertex classes $X_1,\dots,X_\ell$ such that $|X_i|\le cn$ for
all $i=1,\dots,\ell$. Suppose also that $\G$ is a
$(\mathbf{d},\delta_k,\delta,r)$-regular $(k,\ell)$-complex with
vertex classes $V_1,\dots,V_\ell$, all of size $n$, which respects
the partition of~$\Hy$. Then for every vertex $h$ of $\Hy$ we have
that
$$
|\Hy|_{\G} \geq (1 - \alpha) n \left(\prod_{i=2}^k
d_i^{e_i(\Hy)-e_i(\Hy_h)}\right) |\Hy_h|_{\G},
$$
where $\Hy_h$ denotes the induced subcomplex of~$\Hy$ obtained by
removing $h$. In particular, $\G$ contains at least
$((1-\alpha)n)^{|\Hy|}  \prod_{i=2}^k d_i^{e_i(\Hy)}$
labelled partition-respecting copies of $\Hy$.
\end{theorem}
As discussed in the next subsection, Theorem~\ref{emblemma} is a generalization of the
hypergraph counting lemma (which counts subcomplexes~$\Hy$ of bounded order) to
subcomplexes~$\Hy$ of bounded degree and linear order.
Note that the bound relating $|\Hy|_\G$ to $|\Hy_h|_\G$ in
Theorem~\ref{emblemma} is close to what one would get with high
probability if~$\G$ were a random complex\footnote{That is, $\G_2$ is an $\ell$-partite
random graph with density $d_2$, each triangle of $\G_2$ is an edge of $\G_3$ with probability $d_3$ etc.}.
This also shows that
the bound is close to best possible. Theorem~\ref{emblemma} will
be proved in Section~\ref{sec:embed}. In the proof we will need
two lemmas on embeddings of complexes of bounded order, which are
stated in the next subsection.

Recall that if the maximum degree of a $k$-uniform hypergraph~$\Hy$ is at
most $\Delta$ then the maximum degree of the corresponding $k$-complex
is at most~$\Delta2^k$. So it
is easy to see that Theorem~\ref{emblemma} does indeed imply Theorem~\ref{embhgph}.

\subsection{Counting lemma and extension lemma}
We will need a variant (Lemma~\ref{countinglemma}) of the
counting lemma for $k$-unifom hypergraphs due to R\"odl and
Schacht~\cite[Thm~9]{roedlschacht2}. (A similar result was proved
earlier by Gowers~\cite{Gowers} as well as Nagle, R\"odl and Schacht~\cite{Count}.)
It states that if~$|\Hy|$
is bounded and~$\G$ is suitably regular, then the number of copies
of~$\Hy$ in~$\G$ is as large as one would expect if~$\G$ were
random. The main difference to the result in~\cite{roedlschacht2}
is that Lemma~\ref{countinglemma} counts copies of $k$-complexes~$\Hy$ instead of
copies of $k$-uniform hypergraphs~$\Hy$ and also includes an upper bound on the number
of these copies. We will derive Lemma~\ref{countinglemma} from the result
in~\cite{roedlschacht2} in Section~\ref{sec:CL}.

\begin{lemma}[Counting lemma]\label{countinglemma}
Let $k,\ell,r,t,n_0$ be positive integers and let $\eps,d_2,\ldots,
d_k,\delta,\delta_k$ be positive constants such that $1/d_i \in\mathbb{N}$ and
$$1/n_0\ll 1/r, \delta \ll\min\{\delta_k,d_2,\ldots,d_{k-1}\}\le\delta_k\ll
\eps, d_k,1/\ell,1/t.$$%
Then the following holds for all integers $n\ge n_0$. Suppose that
$\Hy$ is a $(k,\ell)$-complex on $t$ vertices with vertex classes
$X_1,\dots,X_\ell$. Suppose also that $\G$ is a
$(\mathbf{d},\delta_k,\delta,r)$-regular $(k,\ell)$-complex with
vertex classes $V_1,\dots,V_\ell$, all of size $n$, which respects
the partition of~$\Hy$. Then
$$|\Hy|_\G=(1\pm\eps) n^{t} \prod_{i=2}^k d_i^{e_i(\Hy)}.$$
\end{lemma}

The main difference between the counting lemma and Theorem~\ref{emblemma}
is that the counting lemma only allows for complexes~$\Hy$
of bounded order. We will apply the counting lemma to embed complexes
of order~$\le f(\Delta,k)$ for some appropriate function $f$.
Note that the upper and lower bounds of the counting lemma imply
Theorem~\ref{emblemma} for the case when~$|\Hy|$ is bounded. A formal
proof of this (which settles the base case for the induction in the proof of
Theorem~\ref{emblemma}) can be found at the beginning of Section~\ref{sec:embed}.

In the induction step of the proof of Theorem~\ref{emblemma} we
will also need the following extension lemma, which states that if
$\Hy'$ is a complex of bounded order, $\Hy\subseteq \Hy'$ is an
induced subcomplex and~$\G$ is suitably regular, then almost all
copies of~$\Hy$ in~$\G$ can be extended to about the `right'
number of copies of~$\Hy'$, where the `right' number is the number
one would expect if~$\G$ were random. We will derive
Lemma~\ref{extensions_count} from Lemma~\ref{countinglemma} in
Section~\ref{sec:ext}.

\begin{lemma}[Extension lemma]\label{extensions_count}
Let $k,\ell,r,t,t',n_0$ be positive integers, where $t<t'$, and let
$\beta,d_2,\ldots,d_k,\delta,\delta_k$ be
positive constants such that $1/d_i \in\mathbb{N}$ and
$$1/n_0\ll 1/r, \delta \ll\min\{\delta_k,d_2,\ldots,d_{k-1}\}\le\delta_k\ll
\beta,d_k,1/\ell,1/t'.$$%
Then the following holds for all integers $n\ge n_0$. Suppose that
$\Hy'$ is a $(k,\ell)$-complex on $t'$ vertices with vertex classes
$X_1,\dots,X_\ell$ and let~$\Hy$ be an induced subcomplex of~$\Hy'$
on~$t$ vertices. Suppose also that $\G$ is a
$(\mathbf{d},\delta_k,\delta,r)$-regular $(k,\ell)$-complex with
vertex classes $V_1,\dots,V_\ell$, all of size $n$, which respects
the partition of $\Hy'$. Then all but at most $\beta |\Hy|_{\G}$
labelled partition-respecting copies of~$\Hy$ in~$\G$ are extendible
to
$$(1 \pm \beta) n^{t'-t} \prod_{i=2}^k d_i^{e_i(\Hy')-e_i(\Hy)}$$
labelled partition-respecting copies of~$\Hy'$ in~$\G$.
\end{lemma}

As well as these versions of the counting lemma and extension
lemma, we will need to be able to apply versions of these lemmas
to underlying $(k-1)$-complexes. In this case, we have that
the regularity constant~$\delta$ is much smaller than all the
densities $d_2,\ldots,d_{k-1}$, but on the other hand we have no~$r$
in the highest level and thus we cannot apply Lemmas~\ref{countinglemma}
and~\ref{extensions_count}. So instead of Lemma~\ref{countinglemma} we
will use the following variant of a result of Kohayakawa, R\"{o}dl
and Skokan~\cite[Cor.~6.11]{KRS}.

\begin{lemma}[Dense counting lemma]\label{densecounting}
Let $k,\ell,t,n_0$ be positive integers and let
$\eps,d_2,\ldots,d_{k-1},\delta$ be positive constants such that
$$ 1/n_0 \ll \delta \ll \eps \ll d_2,\ldots,d_{k-1},1/\ell,1/t.
$$
Then the following holds for all integers $n\ge n_0$. Suppose that
$\Hy$ is a $(k-1,\ell)$-complex on~$t$ vertices with vertex classes
$X_1,\dots,X_\ell$. Suppose also that $\G$ is a
$(d_{k-1},\ldots,d_2,\delta,\delta,1)$-regular $(k-1,\ell)$-complex with
vertex classes $V_1,\dots,V_\ell$, all of size $n$, which respects
the partition of~$\Hy$. Then
$$|\Hy|_\G = (1 \pm \eps) n^t \prod_{i=2}^{k-1}d_i^{e_i(\Hy)}.$$
\end{lemma}

In Section~\ref{sec:CL} we will show how Lemma~\ref{densecounting}
can be deduced from the result in~\cite{KRS}. The following
dense version of the extension lemma can be deduced from the dense
counting lemma (see Section~\ref{sec:ext}).

\begin{lemma}[Dense extension lemma]\label{denseextension}
Let $k,\ell,t,t',n_0$ be positive integers and let
$\beta,d_2,\ldots,d_{k-1},\delta$ be positive constants such that
$$ 1/n_0 \ll \delta \ll \beta \ll d_2,\ldots,d_{k-1},1/\ell,1/t'.
$$
Then the following holds for all integers $n\ge n_0$. Suppose that
$\Hy'$ is a $(k-1,\ell)$-complex on~$t'$ vertices with vertex classes
$X_1,\dots,X_\ell$ and let~$\Hy$ be an induced subcomplex of~$\Hy'$
on~$t$ vertices. Suppose also that $\G$ is a
$(d_{k-1},\ldots,d_2,\delta,\delta,1)$-regular $(k-1,\ell)$-complex with
vertex classes $V_1,\dots,V_\ell$, all of size $n$, which respects
the partition of~$\Hy'$. Then all but at most $\beta |\Hy|_\G$ labelled
partition-respecting copies of~$\Hy$ in~$\G$ can be extended
into
$$ (1\pm \beta) n^{|\Hy'|-|\Hy|}  \prod_{i=2}^{k-1}d_i^{e_i(\Hy')-e_i(\Hy)}$$%
labelled partition-respecting copies of~$\Hy'$ in~$\G$.
\end{lemma}
An overview of how all these lemmas are used in the proof of Theorem~\ref{Ramseythm} is shown
in Figure~1.
\begin{figure}[htb!]  \footnotesize
\centering
\psfrag{1}[][]{Theorem~\ref{cntlemRS}~\cite{roedlschacht2}} \psfrag{2}[][]{Lemma~\ref{DifDensitiesCL}}
\psfrag{3}[][]{Lemma~\ref{countinglemma}} \psfrag{4}[][]{Lemma~\ref{extensions_count}}
\psfrag{5}[][]{Lemma~\ref{densecounting}} \psfrag{6}[][]{Lemma~\ref{denseextension}}
\psfrag{7}[][]{Lemma \ref{KRSthm}~\cite{KRS}}
\psfrag{8}[][]{Theorem~\ref{emblemma}}
\psfrag{9}[][]{Theorem~\ref{reglemma}~\cite{roedlschacht}}  \psfrag{10}[][]{Theorem~\ref{embhgph}}
\psfrag{11}[][]{Proposition~\ref{badktuples}} \psfrag{12}[][]{Theorem~\ref{Ramseythm}}
\psfrag{13}[][]{Lemma~\ref{Slice}}
\includegraphics[scale=0.32]{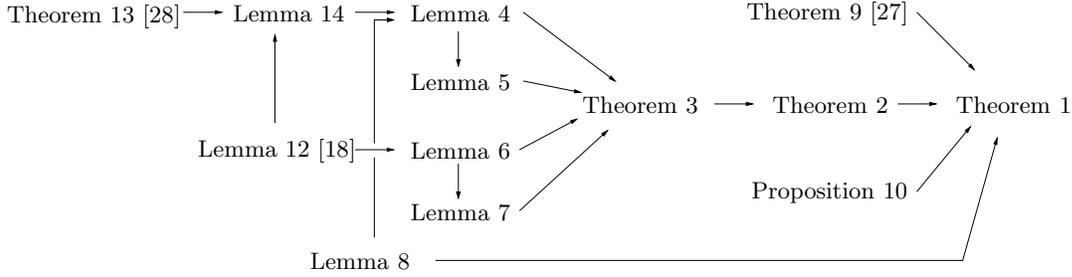}
\caption{Proof of Theorem~\ref{Ramseythm} - Flowchart}
\end{figure}

Another auxiliary result that we will use in the proof of
Lemma~\ref{countinglemma} as well as in the proof of
Theorem~\ref{Ramseythm} is the slicing lemma. Roughly speaking,
this says that in a regular complex $\G$, we can partition the
edge set~$E_j(\G)$ of the $j$th level into an arbitrary number of
parts so that each part is still regular with respect
to~$\G_{j-1}$ with the appropriate density, at the expense of a
larger regularity constant. This can be proved using a simple
application of a Chernoff bound.
\begin{lemma} [Slicing lemma~\cite{roedlschacht}]\label{Slice}
Let $j\geq 2$ and $s_0,r\geq 1$ be integers and let $\delta_0,
d_0$ and $p_0$ be positive real numbers. Then there is an integer
$n_0 = n_0 (j,s_0,r,\delta_0, d_0, p_0)$ such that the following
holds. Let $n \geq n_0$ and let $\G_j$ be a $j$-partite
$j$-uniform hypergraph with vertex classes $V_1,\ldots, V_j$ which
all have size $n$. Also let $\G_{j-1}$ be a $j$-partite
$(j-1)$-uniform hypergraph with the same vertex classes and assume
that each $j$-set of vertices that spans a hyperedge in $\G_j$
also spans a $K_{j}^{(j-1)}$ in $\G_{j-1}$. Suppose that
\begin{enumerate}
\item [1.] $|K_{j-1}^{(j)}(\G_{j-1})| > n^j / \ln n$ and
\item [2.] $\G_j$ is $(d,\delta, r)$-regular with respect to $\G_{j-1}$, where $d\geq d_0 \geq 2\delta \geq 2\delta_0$.
\end{enumerate}
Then for any positive integer $s\leq s_0$ and all positive reals $p_1,\ldots, p_s\ge p_0$ with
$\sum_{i=1}^s p_i \leq 1$ there exists a partition of $E(\G_{j})$ into
$s+1$ parts $E^{(0)}(\G_j), E^{(1)}(\G_j),\ldots, E^{(s)} (\G_j)$ such that
if $\G_{j} (i)$ denotes the spanning subhypergraph of $\G_j$ whose edge set is $E^{(i)}(\G_j)$, then
$\G_{j}(i)$ is $(p_id,3\delta, r)$-regular with respect to $\G_{j-1}$
for every $i=1,\ldots, s$.

Moreover, $\G_{j}(0)$ is
$\left(\left(1-\sum_{i=1}^sp_i\right)d,3\delta, r\right)$-regular
with respect to $\G_{j-1}$ and $E^{(0)}(\G_j) = \emptyset$ if
$\sum_{i=1}^s p_i =1$.
\end{lemma}
%%%%%%%%%%%%%%%%%%%%%%%%%%%%%%%%%%%%%%%%%%%%%%%%%%%%%%%%%%%%%%%%%55
%%%%%%%%%%%%%%%%%%%%%%%%%%%%%%%%%%%%%%%%%%%%%%%%%%%%%%%%%%%%%%%%%%%

\section{Proof of the embedding theorem for complexes (Theorem~\ref{emblemma})}\label{sec:embed}

Throughout the rest of the paper, whenever we talk about a copy of
a complex~$\Hy$ in~$\G$ we mean that this copy is labelled and
partition-respecting, without mentioning this explicitly. We prove
Theorem~\ref{emblemma} by induction on $|\Hy|$. \cite{CNKO}
contains a sketch of the argument for the graph case which gives a
good idea of the proof. We first suppose that the connected
component of~$\Hy$ which contains the vertex~$h$ has order less
than~$\Delta^5$. In this case we will use the counting lemma to
prove the embedding theorem. So let~$\C$ be the component of~$\Hy$
containing~$h$, and let $\D:= \Hy - \C$. Also, let $\C_h:= \C -
h$. We may assume that both $\C_h$ and $\D$ are non-empty. (If
$\D$ is empty then the result follows from
Lemma~\ref{extensions_count}, and if $\C_h$ is empty then $h$ is
an isolated vertex and the result is trivial.) Note that a copy
of~$\Hy$ consists of disjoint copies of~$\C$ and~$\D$,
while~$\Hy_h$ consists of disjoint copies of~$\C_h$ and~$\D$.
Copies of these complexes in~$\G$ will be denoted by~$C$, $D$
and~$C_h$.

Choose a new constant~$\beta$ such that $c,\delta_k\ll\beta\ll\alpha$.
Now note that $|\Hy|_\G = \sum_{D\subseteq \G} |\C|_{\G-D}$, and by
applying the upper and lower bounds of the counting lemma to copies of~$\C$
in~$\G$ and~$\G-D$ respectively, we obtain%
     \COMMENT{Do we need some extra error in case~$\G-D$ is not
regular enough? We do, but it doesn't make things clearer if we include it
and everything is fine if we apply the counting lemma with a
smaller~$\beta$.}
$|\C|_{\G-D} \geq
\frac{(1-c)^{\Delta^5}(1-\beta)}{(1+\beta)}|\C|_{\G}\ge (1-3\beta)|\C|_{\G}$.
So
\begin{equation}\label{basehcount}
|\Hy|_\G \geq \sum_{D\subseteq \G} (1-3\beta)|\C|_{\G}=
(1-3\beta)|\C|_{\G} |\D|_\G.
\end{equation}
On the other hand, by a similar argument using the upper and
lower bounds from the counting lemma in~$\G$ for~$\C_h$ and~$\C$
respectively,
\begin{equation}\label{basehhcount}
|\Hy_h|_\G \leq |\C_h|_\G |\D|_\G \leq
\frac{1+\beta}{1-\beta}\frac{|\C|_\G |\D|_\G}{n\prod_{i=2}^k
d_i^{e_i(\C)-e_i(\C_h)}}.
\end{equation}
Combining (\ref{basehcount}) and (\ref{basehhcount}) gives the
desired result.

Thus we may assume that the component of $\Hy$ containing~$h$
has order at least~$\Delta^5$. This deals with the base case
of the inductive argument, and it also means that the fourth
neighbourhood of~$h$ in~$\Hy$ will be non-empty, which will be
convenient later on in the proof as we will only be counting complexes which
are non-empty.

We pick new constants $\eps_k$ and $\eps_{k-1}$ satisfying the following
hierarchies:%
$$ \delta \ll \eps_{k-1} \ll d_2,d_3,\ldots,d_k, 1/\Delta,$$
$$ c,\delta_k,\eps_{k-1} \ll \eps_k \ll \alpha.
$$
Let~$\N_h$ be the subcomplex of~$\Hy$ induced by the neighbours of~$h$,
and let~$\B$ be the subcomplex of~$\Hy$ induced by~$h$ and
the neighbours of~$h$. Then any copy of~$\Hy$ in~$\G$ extending
a copy~$N_h$ of~$\N_h$ can be obtained by first extending~$N_h$ into
a copy of~$\Hy_h$ and then extending~$N_h$ into a copy of~$\B$,
where the vertex chosen for~$h$ has to be distinct from all the
vertices chosen for~$\Hy_h$.

We now introduce some more notation. Given $k$-complexes $\Hy'
\subseteq \Hy''$ such that~$\Hy'$ is induced, and a copy~$H'$
of~$\Hy'$ in~$\G$, we define $\hptohpp$ to be the number of ways
in which~$H'$ can be extended to a copy of $\Hy''$ in $\G$. We
also define
$$\exphptohpp := n^{|\Hy''|-|\Hy'|}\prod_{i=2}^k d_i^{e_i(\Hy'')-e_i(\Hy')}.
$$
Thus $\exphptohpp$ is roughly the expected number of ways $H'$
could be extended to a copy of $\Hy''$ if $\G$ were a random
complex.

We define a copy~$N_h$ of~$\N_h$ to be \emph{typical} if it
has about the correct number of extensions into~$\B$, i.e. if $\nhtob
= (1 \pm \eps_k)\expnhtob$. An application of the extension lemma
(Lemma~\ref{extensions_count})
shows that at most $\eps_k |\N_h|_\G$ copies of~$\N_h$ in~$\G$ are not
typical. We denote the set of typical copies of~$\N_h$ by~$\typ$,
and the set of all atypical copies by~$\atyp$.

Now observe that if all of the copies of $\N_h$ were typical, the
proof would be complete, since then
\begin{align} \label{eq:BasicTyp}
|\Hy|_\G & \geq \sum_{N_h\subseteq \G} \nhtohh (\nhtob - cn)
 \geq \left((1-\eps_k)\expnhtob - cn\right)\sum_{N_h \subseteq \G} \nhtohh \nonumber\\
& \geq (1-\alpha) \expnhtob |\Hy_h|_\G
= (1-\alpha) \exphhtoh |\Hy_h|_\G.
\end{align}
The third inequality follows since $c \ll \alpha,d_2,\ldots,d_k$,
and $\eps_k \ll \alpha$.

However, we also need to take account of the atypical copies of~$\N_h$.
The proportion of these is about~$\eps_k$, which may be
larger than some~$d_i$. It will turn out that this is too large for
our purposes, and so we will need to consider the atypical copies
more carefully.

We define, instead of $\hptohpp$, the
expression $\hptohppmin$, where $\Hy'\subseteq \Hy''$ are induced subcomplexes
of~$\Hy$ and~$H'$ is a copy of~$\Hy'$ in~$\G$.
We consider the underlying $(k-1)$-complexes in each case, and define
$\hptohppmin$ to be the number of ways in which the underlying
$(k-1)$-complex of~$H'$ can be extended to the underlying
$(k-1)$-complex of~$\Hy''$ within (the underlying $(k-1)$-complex
of)~$\G$. Clearly $\hptohppmin \geq \hptohpp$. We also define
$$\exphptohppmin := n^{|\Hy''|-|\Hy'|} \prod_{i=2}^{k-1}
d_i^{e_i(\Hy'')-e_i(\Hy')}.
$$
Thus $\exphptohppmin$ is roughly the expected value of
$\hptohppmin$ if~$\G$ were a random complex. Also,
$$
\exphptohppmin
= \exphptohpp/d_k^{e_k(\Hy'')-e_k(\Hy')}\geq \exphptohpp.
$$

We define~$\N_h^*$ to be the subcomplex of~$\Hy$ induced by
the vertices at distance~$3$ from~$h$. We also define~$\F$ to be
the subcomplex of~$\Hy$ induced by the vertices at distance
$1,2$ or $3$ from~$h$, i.e. the subcomplex induced by~$\N_h$,
$\N_h^*$ and the vertices in between (see Figure~2).
\begin{figure}[htb!]  \footnotesize
\centering
\psfrag{1}[][]{$h$} \psfrag{2}[][]{$\ \cn$}
\psfrag{3}[][]{$\ \cn^*$}
\psfrag{5}[][]{$\Hy_h^{*}$}
\psfrag{4}[][]{$\F$} \psfrag{6}[][]{$\Hy^-_h$}
\psfrag{7}[][]{$\Hy_h$}
\psfrag{8}[][]{$\Hy'_h$}
\psfrag{9}[][]{$\B$}
\includegraphics[scale=0.35]{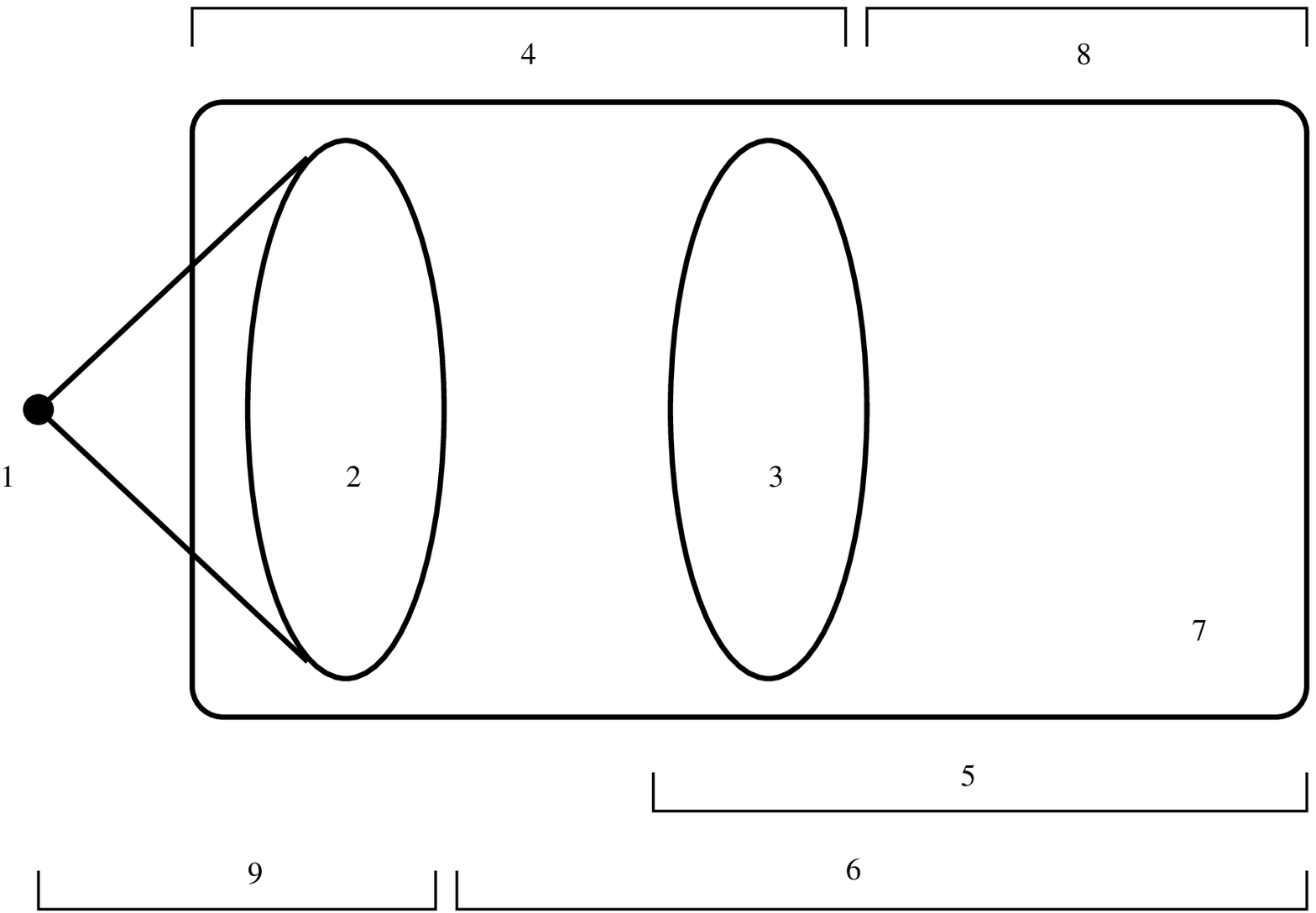}
\caption{The complex $\Hy$}
\end{figure}
%Given disjoint copies~$N_h$ of~$\N_h$ and~$N_h^*$ of~$\N_h^*$ we
%define $|N_h,N_h^* \rightarrow \F|:= |N_h\cup N_h^*\rightarrow
%\F|$ to be the number of ways $N_h\cup N_h^*$ can be extended into
%a copy of~$\F$. We set $\overline{|\N_h,\N_h^*\rightarrow
%\F|}:= \overline{|\N_h\cup\N_h^*\rightarrow \F|}$. The
%expressions $\nhnhstofmin$ and $\expnhnhstofmin$ are defined
%similarly.

Given copies~$N_h$ of~$\N_h$ and~$N_h^*$ of~$\N_h^*$, we say that
the pair $N_h,N_h^*$ is \emph{useful} if~$N_h$ and~$N_h^*$ are
disjoint and if the pair has about the expected number of extensions
into copies of~$\F$ as $(k-1)$-complexes, i.e. if
$$\nhnhstofmin =(1 \pm \eps_{k-1})\expnhnhstofmin. $$

We use Lemmas~\ref{countinglemma}, \ref{densecounting}
and~\ref{denseextension} applied to $\N_h \cup \N_h^*$ to show
that at most $\sqrt{\eps_{k-1}}|\N_h|_\G|\N_h^*|_\G$ disjoint
pairs $N_h,N_h^*$ are not useful. Let $|\N_h \cup
\N_h^*|_\G^{(k-1)}$ denote the number of copies of the underlying
$(k-1)$-complex of~$\N_h \cup \N_h^*$ in~$\G$. Then
Lemmas~\ref{countinglemma} and~\ref{densecounting} together imply
that $|\N_h \cup \N_h^*|_\G^{(k-1)}\le (1+2\eps_k) |\N_h \cup
\N_h^*|_\G /d_k^{e_k(\N_h\cup \N_h^*)}$. Moreover, the dense
extension lemma (Lemma~\ref{denseextension}) shows that all but at
most $\eps_{k-1}|\N_h \cup \N_h^*|_\G^{(k-1)}$ copies of the
underlying $(k-1)$-complex of~$\N_h\cup \N_h^*$ in~$\G$ are
useful. Altogether this shows that all but at most
\begin{equation} \label{upperuseful}
\eps_{k-1}(1+2\eps_k) |\N_h\cup \N_h^*|_\G /d_k^{e_k(\N_h\cup \N_h^*)}\le \sqrt{\eps_{k-1}}|\N_h \cup \N_h^*|_\G
\leq \sqrt{\eps_{k-1}}|\N_h |_{\G} | \N_h^*|_\G
\end{equation}
disjoint pairs of copies of~$\N_h$ and $\N_h^*$ are useful. Note
that if we had chosen~$\N_h^*$ to be the subcomplex of $\Hy$
induced by the vertices at distance 2 from $h$ (instead of 3),
then we could not have applied Lemma~\ref{densecounting}, since
$\N_h \cup \N_h^*$ would not be an induced subcomplex of $\F$.
Together with the fact that only comparatively few of the pairs
$N_h,N_h^*$ will intersect, this shows that at most
$2\sqrt{\eps_{k-1}}|\N_h|_\G|\N_h^*|_\G$ pairs $N_h,N_h^*$ are not
useful. Hence at most $\eps_{k-1}^{1/4}|\N_h|_\G$ copies of~$\N_h$
form a non-useful pair together with more than
$2\eps_{k-1}^{1/4}|\N_h^*|_\G$ copies of~$\N_h^*$. We call all
other copies of $\N_h$ \emph{useful} and let~$\usef$ denote the
set of all these copies. Then
\begin{equation} \label{nonusefbound}
|\N_h|_\G -|\usef| \leq  \eps_{k-1}^{1/4} |\N_h|_\G.
\end{equation}
We denote by $\usef^*(N_h)$ the set of all~$N_h^*$ which form
a useful pair together with~$N_h$.

\medskip

\noindent\textbf{Claim.} \emph{Any useful copy $N_h$ of $\N_h$ satisfies
\begin{equation*}
\nhtohh \leq \frac{10}{d_k^{\Delta^3}}
\frac{|\Hy_h|_\G}{|\N_h|_\G}.
\end{equation*}}

\smallskip

\noindent
Note that $\sum_{N_h} \nhtohh = |\Hy_h|_\G$, so $|\Hy_h|_\G / |\N_h|_\G$ is the average
value of $\nhtohh$ over all copies~$N_h$ of~$\N_h$.
Later on, we will apply the claim to show that only a
small fraction of copies of~$\Hy$ contain a useful but atypical
copy of~$\N_h$.

\medskip

\noindent
{\bf Proof of Claim.} Fix a useful copy~$N_h$
of~$\N_h$. Put~$\Hy_h^*:=\Hy_h-(\F-\N_h^*)$. We aim to
extend~$N_h$ to a copy of~$\Hy_h$ by first picking a copy~$N_h^*$
of~$\N_h^*$, then extending this to a copy of~$\Hy_h^*$ and also
extending $N_h\cup N_h^*$ to a copy of~$\F$. We must also make sure
that no vertices are used more than once. However, since we are
only looking for an upper bound on $\nhtohh$, and ignoring this
restriction can only increase the number of extensions we find, we
may ignore this difficulty. Thus
\begin{equation}\label{usefsums}
\nhtohh \leq \sum_{N_h^* \in \usef^*(N_h)} \!\!\!\!\!\!\! \nhnhstof \nhstohhs +
\sum_{N_h^* \notin \usef^*(N_h)} \!\!\!\!\!\!\! \nhnhstof \nhstohhs.
\end{equation}
%We know that only a small fraction ($\eps_{k-1}$) of copies of
%$\N_h^*$ are not useful. We will therefore be able to use the crude
%bound $\nhstohhs \leq n^{|\Hy_h^*|-|\N_h^*|}$.
We bound the two sums separately. To bound the first sum, we need to
bound $\nhnhstof$ in the case when the pair $N_h,N_h^*$ is useful.
But clearly $\nhnhstof \leq \nhnhstofmin$, and
$$
\nhnhstofmin \leq (1+\eps_{k-1})\expnhnhstofmin
=\frac{(1+\eps_{k-1})\expnhnhstof}{d_k^{e_k(\F)-e_k(\N_h)-e_k(\N_h^*)}}
$$
whenever $N_h^* \in \usef^*(N_h)$.
So the first sum in~(\ref{usefsums}) is bounded by
\begin{equation} \label{Hupperbound}
\frac{1+\eps_{k-1}}{d_k^{e_k(\F)-e_k(\N_h)-e_k(\N_h^*)}}
\expnhnhstof|\Hy_h^*|_\G \le \frac{2}{d_k^{\Delta^3}}
\expnhnhstof|\Hy_h^*|_\G.
\end{equation}
To see the bound of $\Delta^3$ on the number of $k$-edges which we used
in the final inequality, note that
$|\F-\N_h-\N_h^*| \le \Delta^2$ and that the number of $k$-edges each of these vertices
lies in is at most $\Delta$.
We now want to express the bound in~(\ref{Hupperbound}) in terms of $|\Hy_h^-|_\G$,
where $\Hy_h^-:=\Hy_h-\N_h$. By the induction hypothesis applied
several times,
\begin{align*}
|\Hy_h^*|_\G & \leq ((1-\alpha)n)^{-(|\Hy_h^-|-|\Hy_h^*|)}
\left(\prod_{i=2}^k d_i^{-(e_i(\Hy_h^-)-e_i(\Hy_h^*))}\right)
|\Hy_h^-|_\G\\
& \le 2 \frac{\prod_{i=2}^k
d_i^{e_i(\Hy_h)-e_i(\Hy_h^-)-e_i(\N_h)}}{\expnhnhstof}|\Hy_h^-|_\G.
\end{align*}
In the last line we used that $e_i(\Hy_h)=e_i(\Hy_h^*)+e_i(\F)- e_i(\N_h^*)$
and $|\F|-|\N_h|-|\N_h^*|=|\Hy_h^-|-|\Hy_h^*|$ (see Figure~1).
We also used that $(1-\alpha)^{-(|\Hy_h^-|-|\Hy_h^*|)} \le 2$.
So we obtain
\begin{align}\label{usef}
\sum_{N_h^* \in \usef^*(N_h)} \nhnhstof \nhstohhs  & \leq
\frac{4 \prod_{i=2}^k
d_i^{e_i(\Hy_h)-e_i(\Hy_h^-)-e_i(\N_h)} }
{d_k^{\Delta^3}} |\Hy_h^-|_\G .
\end{align}
To bound the second sum in~(\ref{usefsums}), we define $\Hy_h':=
\Hy_h^* - \N_h^*$, and observe that trivially any copy
$N_h^*$ of $\N_h^*$ satisfies $\nhstohhs \leq |\Hy_h'|_\G$. Note that $\Hy_h'$ is nonempty.
On the other hand, by the induction hypothesis
applied several times,%
\begin{align*}
|\Hy_h'|_\G & \leq
((1-\alpha)n)^{|\Hy_h'|-|\Hy_h^-|}\left(\prod_{i=2}^k
d_i^{e_i(\Hy_h')-e_i(\Hy_h^-)}
\right)|\Hy_h^-|_\G \leq
\frac{2|\Hy_h^-|_\G}{\left(\prod_{i=2}^kd_i\right)^{2\Delta^4}
n^{|\Hy_h^-|-|\Hy_h'|}}.
\end{align*}
Since at most $2\eps_{k-1}^{1/4}|\N_h^*|_\G \le 2\eps_{k-1}^{1/4}n^{|\N_h^*|}$
copies of~$\N_h^*$ do not lie
in~$\usef(N_h)$, the second sum in~(\ref{usefsums}) is bounded by%
\begin{eqnarray}\label{non-usef}
\sum_{N_h^* \notin \usef(N_h)} \!\!\!\!\!\!\! \nhnhstof \nhstohhs & \leq &
2\eps_{k-1}^{1/4} n^{|\N_h^*|}n^{|\F|-|\N_h|-|\N_h^*|}
\frac{2|\Hy_h^-|_\G}{\left(\prod_{i=2}^k d_i\right)^{2\Delta^4}
n^{|\Hy_h^-|-|\Hy_h'|}} \nonumber \\
& = & 2\eps_{k-1}^{1/4}
\frac{2|\Hy_h^-|_\G}{\left(\prod_{i=2}^k d_i\right)^{2\Delta^4}
} \nonumber \\
& \leq & \left(\prod_{i=2}^k
d_i^{e_i(\Hy_h)-e_i(\Hy_h^-)-e_i(\N_h)}\right)|\Hy_h^-|_\G.
\end{eqnarray}
The last inequality follows since $\eps_{k-1} \ll
d_2,d_3,\ldots,d_k,1/\Delta$. Substituting (\ref{usef}) and
(\ref{non-usef}) into (\ref{usefsums}) we obtain
\begin{eqnarray}\label{usefexthhmin}
\nhtohh & \leq &
\left(1+\frac{4}{d_k^{\Delta^3}}\right)\left(\prod_{i=2}^k
d_i^{e_i(\Hy_h)-e_i(\Hy_h^-)-e_i(\N_h)}\right)|\Hy_h^-|_\G \nonumber \\
& \leq & \frac{5\left(\prod_{i=2}^k
d_i^{e_i(\Hy_h)-e_i(\Hy_h^-)-e_i(\N_h)}\right)}{d_k^{\Delta^3}}|\Hy_h^-|_\G.
\end{eqnarray}
It now remains only to relate $|\Hy_h^-|_\G$ to
$|\Hy_h|_\G/|\N_h|_\G$. Once again we apply the induction hypothesis
several times to obtain%
\begin{equation*}
|\Hy_h|_\G \geq
((1-\alpha)n)^{|\Hy_h|-|\Hy_h^-|}\prod_{i=2}^k
d_i^{e_i(\Hy_h)-e_i(\Hy_h^-)}|\Hy_h^-|_\G.
\end{equation*}
On the other hand, the counting lemma implies that
$|\N_h|_\G \leq (1+\alpha)\left(\prod_{i=2}^k d_i^{e_i(\N_h)}\right)n^{|\N_h|}$.
Putting these two bounds together, we obtain%
\begin{eqnarray}\label{hhtohhminnh}
\frac{|\Hy_h|_\G}{|\N_h|_\G} & \geq &
\frac{((1-\alpha)n)^{|\Hy_h|-|\Hy_h^-|}\left(\prod_{i=2}^k
d_i^{e_i(\Hy_h)-e_i(\Hy_h^-)}\right)|\Hy_h^-|}{(1+\alpha)\left(\prod_{i=2}^k
d_i^{e_i(\N_h)}\right)n^{|\N_h|}} \nonumber \\
& \geq & \frac{1}{2} \left(\prod_{i=2}^k
d_i^{e_i(\Hy_h)-e_i(\Hy_h^-)-e_i(\N_h)}\right)|\Hy_h^-|_\G.
\end{eqnarray}
Together with (\ref{usefexthhmin}), this shows that
\begin{align*}
\nhtohh & \leq \frac{5 \cdot 2}{d_k^{\Delta^3}}\frac{|\Hy_h|_\G}{|\N_h|_\G},
\end{align*}
which completes the proof of the claim.
\endproof
\medskip

\noindent
Using the claim we now go on to prove the induction step.
Given a copy~$H_h$ of~$\Hy_h$, we denote by~$N_h(H_h)$ the induced
copy of~$\N_h$. We have
\begin{eqnarray}\label{prelimhcount}
|\Hy|_\G & = & \sum_{H_h\subseteq \G} \hhtoh
\geq \sum_{H_h\subseteq \G} (|N_h(H_h) \rightarrow \B| -cn)\nonumber \\
& = & \sum_{N_h\subseteq \G} \nhtohh \nhtob - cn|\Hy_h|_\G \nonumber \\
& \geq & (1-\eps_k)\expnhtob \left(\sum_{N_h\subseteq \G} \nhtohh -
\sum_{N_h \notin \typ} \nhtohh \right) - cn|\Hy_h|_\G.
\end{eqnarray}
We want to show that the term in this expression which comes from
the atypical copies of~$\N_h$ does not affect the calculations too
much, and so we aim to bound the contribution from atypical copies
of~$\N_h$. We have
\begin{equation}\label{atypsum}
\sum_{N_h \notin \typ} \nhtohh = \sum_{N_h \notin \typ, N_h \in
\usef} \nhtohh + \sum_{N_h \notin \typ, N_h \notin
\usef} \nhtohh .
\end{equation}
Now the claim implies that we can bound the first sum
in~(\ref{atypsum}) by
\begin{eqnarray}\label{atypusef}
\sum_{N_h \notin \typ, N_h \in \usef} \nhtohh & \leq & \sum_{N_h
\notin \typ, N_h \in \usef} \frac{10}{d_k^{\Delta^3}}
\frac{|\Hy_h|_\G}{|\N_h|_\G}
\leq |\atyp|
\frac{10}{d_k^{\Delta^3}} \frac{|\Hy_h|_\G}{|\N_h|_\G}  \leq
\sqrt{\eps_k}|\Hy_h|_\G.
\end{eqnarray}
Meanwhile we can also bound the second sum by
\begin{eqnarray}\label{atypnonusef}
\sum_{N_h \notin \typ, N_h \notin
\usef} \nhtohh & \leq & \sum_{N_h \notin \typ, N_h
\notin \usef} |\Hy_h^-|_\G  \nonumber \\
& \stackrel{(\ref{hhtohhminnh})}{\leq} & (|\N_h|_\G-|\usef|)
\frac{2}{\prod_{i=1}^kd_i^{\Delta^2}}
\frac{|\Hy_h|_\G}{|\N_h|_\G} \nonumber \\
& \stackrel{(\ref{nonusefbound})}{\leq} & \eps_{k-1}^{1/5} |\Hy_h|_\G.
\end{eqnarray}
Combining (\ref{atypsum}), (\ref{atypusef}) and
(\ref{atypnonusef}), we have%
$$\sum_{N_h \notin \typ} \nhtohh \leq 2 \sqrt{\eps_k} |\Hy_h|_\G
$$
and combining this with (\ref{prelimhcount}), we obtain
\begin{align*}
|\Hy|_\G & \geq (1-\eps_k) \expnhtob \left( |\Hy_h|_\G -
2\sqrt{\eps_k}
|\Hy_h|_\G \right)- cn|\Hy_h|_\G\\
& = (1-\eps_k)n \left(\prod_{i=2}^k d_i^{e_i(\B)-e_i(\N_h)}
\right) (1-2\sqrt{\eps_k})|\Hy_h|_\G - cn|\Hy_h|_\G\\
& \geq (1-\alpha)n\left( \prod_{i=2}^k d_i^{e_i(\Hy)-e_i(\Hy_h)}
\right)|\Hy_h|_\G,
\end{align*}
as required. This completes the proof of Theorem~\ref{emblemma}.

\section{The regularity lemma for $k$-uniform hypergraphs} \label{regularity}
\subsection{Preliminary definitions and statement}

In this section we state the version of the regularity lemma for $k$-uniform
hypergraphs due to R\"odl and Schacht~\cite{roedlschacht}, which we use
in the proof of Theorem~\ref{Ramseythm} in the next section.
To prepare for this we will first need some notation. We follow~\cite{roedlschacht}.
Given a finite set~$V$ of vertices, we will define a family $\Part =
\{\Part^{(1)},\ldots, \Part^{(k-1)}\}$ where each $\Part^{(j)}$ is a partition
of certain $j$-subsets of~$V$. These partitions will satisfy  properties which
we will describe below. We denote by~$[V]^j$ the set of all $j$-subsets of~$V$.
Suppose that we are given a partition $\Part^{(1)}=\{V_1,\ldots, V_{|\Part^{(1)}|} \}$
of $[V]^1 = V$. We will call the~$V_i$ \emph{clusters}.
We denote by $\cross_j=\cross_j(\Part^{(1)})$ the set
of all those $j$-subsets of~$V$ that meet each part of~$\Part^{(1)}$ in at most~1
element. Each~$\Part^{(j)}$ will be a partition of~$\cross_j$. Moreover,
any two $j$-sets that belong to the same part of~$\Part^{(j)}$ will meet the
same~$j$ clusters. This means that each part of~${\Part}^{(j)}$ can be viewed as
a $j$-partite $j$-uniform hypergraph whose vertex classes are these clusters.
In particular, the parts of~$\Part^{(2)}$ can be thought of as
bipartite subgraphs between two of the clusters. Moreover, for each part~$A$
of~$\Part^{(3)}$ there will be~3 clusters and~3 bipartite graphs belonging
to~$\Part^{(2)}$ between these clusters such that all the 3-sets in~$A$
form triangles in the union of these~3 bipartite graphs.

More generally, suppose that we have already defined partitions
$\Part^{(1)},\ldots, \Part^{(j-1)}$ and are about to define~$\Part^{(j)}$.
Given $i<j$ and $I \in \cross_i$, we let~$P^{(i)}(I)$
denote the part of~$\Part^{(i)}$ the set~$I$ belongs to.
Given $J\in \cross_j$, the \emph{polyad~$\hat{P}^{(j-1)}(J)$ of~$J$} is
defined by
$$ \hat{P}^{(j-1)} (J) := \bigcup \{P^{(j-1)}(I) : \; I \in [J]^{j-1} \}.
$$
Thus~$\hat{P}^{(j-1)}(J)$ is the unique collection of~$j$ parts of~$\Part^{(j-1)}$
in which~$J$ spans a copy of the complete~$(j-1)$-uniform hypergraph~$K_j^{(j-1)}$
on~$j$ vertices. Moreover, note that $\hat{P}^{(j-1)}(J)$ can be viewed as
a $j$-partite $(j-1)$-uniform hypergraph whose vertex classes are the~$j$ clusters
containing the vertices of~$J$. We set
$$ \hat{\Part}^{(j-1)} := \{ \hat{P}^{(j-1)}(J): \; J \in \cross_j \}.
$$
Note that the polyads~$\hat{P}^{(j-1)}(J)$ and~$\hat{P}^{(j-1)}(J')$ need not be
distinct for different $J,J'\in[V]^j$. However, if these polyads are distinct then
$\K_j(\hat{P}^{(j-1)}(J))\cap\K_j(\hat{P}^{(j-1)}(J'))=\emptyset$.
(Recall that $\K_j(\hat{P}^{(j-1)}(J))$ is the set of all $j$-sets of vertices
which form a~$K_j^{(j-1)}$ in~$\hat{P}^{(j-1)}(J)$. So in particular,
$\K_j(\hat{P}^{(j-1)}(J))$ contains~$J$.)
This implies that $\{\K_j(\hat{P}^{(j-1)}): \hat{P}^{(j-1)}\in \hat{\Part}^{(j-1)} \}$
is a partition of $\cross_j$. The property of~$\Part^{(j)}$ which we require
is that it refines $\{\K_j(\hat{P}^{(j-1)}):\; \hat{P}^{(j-1)}\in \hat{\Part}^{(j-1)} \}$,
i.e.~each part of $\Part^{(j)}$ has to be contained in some~$\K_j(\hat{P}^{(j-1)})$.

We also need a notion which generalizes that of a polyad: given $J\in \cross_j$ and $i<j$ we set
$$ \hat{P}^{(i)} (J) := \bigcup \{P^{(i)}(I) : \; I \in [J]^{i} \}.
$$
Then the properties of our partitions imply that
$\bigcup_{i=1}^{j-1}\hat{P}^{(i)}(J)$ is a $(j-1,j)$-complex.

Altogether, given $\av = (a_1,\ldots, a_{k-1}) \in \mathbb{N}^{k-1}$ we say
that $\Part(k-1,\av)=\{\Part^{(1)},\ldots, \Part^{(k-1)} \}$
is a \emph{family of partitions on~$V$} if
\begin{enumerate}
\item $\Part^{(1)}$ is a partition of~$V$ into $a_1$ clusters.
\item For all $j=2,\dots,k-1$, $\Part^{(j)}$ is a partition of $\cross_j$
such that for each part
there is a polyad $\hat{P}^{(j-1)} \in \hat{\Part}^{(j-1)}$ so that the part is
contained in $\K_j(\hat{P}^{(j-1)})$. Moreover, for each polyad $\hat{P}^{(j-1)} \in
\hat{\Part}^{(j-1)}$, the set $\K_j(\hat{P}^{(j-1)})$ is the union of~$a_j$ parts
of~$\Part^{(j)}$.
\end{enumerate}
We say that $\Part =\Part(k-1,\av)$ is \emph{$t$-bounded} if
$a_1,\ldots, a_{k-1}\leq t$.
Suppose that $a_1$ divides $|V|$. Then $\Part = \Part(k-1,\av)$ is called
\emph{$(\eta, \delta,\av)$-equitable} if
\begin{enumerate}
\item $\Part^{(1)}$ is a partition of~$V$ into~$a_1$ clusters of equal size;
\item $|[V]^k \setminus \cross_k| \leq \eta {|V| \choose k}$;
\item for every $K \in \cross_k$, the $(k-1,k)$-complex $\bigcup_{i=1}^{k-1}\hat{P}^{(i)}(K)$
is $({\bf d},\delta,\delta,1)$-regular, where ${\bf d} =(1/a_{k-1},\ldots, 1/a_2)$.
\end{enumerate}
In particular, the second condition implies that $1/a_1$ is small compared to~$\eta$.
%We will interpret this as requiring that
%$1/a_1 \ll \eta$, although the latter is a stronger condition.

Let $\delta_k>0$ and~$r\in\mathbb{N}$. Suppose that~$\G$ is a
$k$-uniform hypergraph on~$V$ and
$\Part = \Part (k-1,\av)$ is a family of partitions on~$V$. Recall that we can view
each polyad~$\hat{P}^{(k-1)}\in \hat{\Part}^{(k-1)}$ as a $(k-1)$-uniform $k$-partite
hypergraph. $\G$ is called \emph{$(\delta_k,r)$-regular with respect
to $\hat{P}^{(k-1)}$} if~$\G$ is $(d,\delta_k,r)$-regular with respect
to $\hat{P}^{(k-1)}$ for some~$d$. We say
that~$\G$ is  \emph{$(\delta_k,r)$-regular with respect to~$\Part$} if
$$\left| \bigcup \{ \K_k(\hat{P}^{(k-1)}): \; \mbox{ $\G$ is not $(\delta_k,r)$-regular with respect
to $\hat{P}^{(k-1)} \in \hat{\Part}^{(k-1)}$}  \} \right| \leq \delta_k |V|^k.
$$
This means that not much more than a $\delta_k$-fraction of the $k$-subsets of~$V$
form a~$K_k^{(k-1)}$ that lies within a polyad with respect to which~$\G$
is not regular.

Now, we are ready to state the regularity lemma, which we are
going to use in the proof of Theorem~\ref{Ramseythm}.

\begin{theorem}[R\"odl and Schacht~\cite{roedlschacht}]\label{reglemma}
Let $k\geq 2$ be a fixed integer. For all positive constants~$\eta$ and~$\delta_k$
and all functions $r: {\mathbb N}^{k-1}\rightarrow \mathbb{N}$ and
$\delta: \mathbb{N}^{k-1} \rightarrow (0,1]$, there are integers~$t$ and~$m_0$ such
that the following holds for all $m\ge m_0$ which are divisible by~$t!$.
Suppose that~$\G$ is a $k$-uniform hypergraph of order~$m$. Then there exists an
$\av\in\mathbb{N}^{k-1}$ and a
family of partitions $\Part = \Part(k-1, \av)$ of the vertex set~$V$ of~$\G$
such that
\begin{enumerate}
\item $\Part$ is $(\eta,\delta (\av),\av)$-equitable and~$t$-bounded and
\item $\G$ is $(\delta_k,r(\av))$-regular with respect to~$\Part$.
\end{enumerate}
\end{theorem}
The advantage of this regularity lemma compared to the one proved earlier by R\"odl
and Skokan~\cite{RSkok} is that it uses only two regularity constants~$\delta$
and~$\delta_k$ instead of~$k-1$ different ones. The regularity constants
$\delta_2,\ldots, \delta_k$ produced by the regularity lemma
in~\cite{RSkok} might satisfy
$\delta_2 \ll 1/a_2 \ll \delta_3  \ll 1/a_3
\ll \cdots \ll 1/a_{k-1} \ll \delta_k$, which would make the proof of the
corresponding embedding theorem more technical in appearance.

Note that the constants in Theorem~\ref{reglemma} can be chosen such that
they satisfy the following hierarchy:
\begin{equation} \label{reghierarchy}
\frac{1}{m_0} \ll \frac{1}{r}=\frac{1}{r(\av)},\delta=\delta(\av)
\ll \min\{\delta_k,\eta,1/a_1,1/a_2,\ldots,1/a_{k-1}\}.
\end{equation}

\subsection{The reduced hypergraph}\label{sec:rhgraph}

In the proof of Theorem~\ref{Ramseythm} that follows in the next
section, we will use the so-called reduced hypergraph. If $\Part =
\{ \Part^{(1)},\ldots, \Part^{(k-1)} \}$ is the partition of the vertex set
of~$\G$ given by the regularity lemma, the \emph{reduced hypergraph}
$\R=\R(\G,\Part)$ is a $k$-uniform hypergraph whose vertices are the
clusters, i.e.~the parts of~$\Part^{(1)}$. To
define the set of hyperedges we need the following notion.
We say that a $k$-tuple of clusters is \emph{fruitful} if~$\G$ is
$(\delta_k,r)$-regular with respect to all but at most a
$\sqrt{\delta_k}$-fraction of all those polyads~$\hat{P}^{(k-1)}$
which are induced on these~$k$
clusters. The set of hyperedges of~$\R$ consists of precisely those
$k$-tuples that are fruitful. In the proof of
Theorem~\ref{Ramseythm}, we shall need an estimate on the number
of these hyperedges. In particular, we need to show that~$\R$ is
very dense. This is conveyed in the following proposition.

\begin{proposition} \label{badktuples}
All but at most $2\sqrt{\delta_k} a_1^k$ of the $k$-tuples of clusters are fruitful.
\end{proposition}
\proof By the dense counting lemma (Lemma~\ref{densecounting})
each polyad in~$\hat{\Part}^{(k-1)}$ contains at least
$$ f(m,\av):=\frac{1}{2} \left( \frac{m}{a_1}\right)^k \prod_{i=2}^{k-1}
\left(\frac{1}{a_i}\right)^{{k \choose i}}
$$
copies of~$K_k^{(k-1)}$. Since~$\G$ is $(\delta_k,r)$-regular with
respect to~$\Part$, the number of polyads in~$\hat{\Part}^{(k-1)}$
with respect to which~$\G$ is not $(\delta_k,r)$-regular is at most
\begin{equation} \label{bound}
\frac{\delta_k m^k}{f(m,\av)} = \frac{2\prod_{i=1}^{k-1}a_i^{{k
\choose i}}}{m^k} \delta_k m^k = 2\delta_k
\prod_{i=1}^{k-1}a_i^{{k \choose i}}.
\end{equation}
We call these polyads \emph{bad}. Now, each $k$-tuple of clusters
induces $\prod_{i=2}^{k-1} a_i^{{k \choose i}}$ polyads in $\hat{\Part}^{(k-1)}$.%
     \COMMENT{If e.g. $k=4$ then to choose a polyad we have to choose
a bipartite graph between each pair of the 4 clusters. So there
are $a_2^{\binom{4}{2}}$ possible choices for these bipartite graphs.
Now for each of the $\binom{4}{3}$ 3-tuples of clusters we have to choose
a part of~$\Part^{3}$ which lies in the set of all triangles spanned
by the chosen bipartite graphs between these 3 clusters. Since the
set of these triangles is partitioned into $a_3$ parts of $\Part^{3}$,
we have $a_3$ choices for such a part. Thus in total for all the 3-tuples
we have $a_3^{\binom{4}{3}}$ choices.}
Thus if there were more than $2\sqrt{\delta_k} a_1^k$ $k$-tuples of
clusters each inducing more than
$\sqrt{\delta_k} \prod_{i=2}^{k-1} a_i^{{k \choose i}}$ bad polyads, the total
number of bad polyads would exceed the bound given
in~(\ref{bound}), yielding a contradiction.
\endproof

\section{Proof of Theorem~\ref{Ramseythm}}\label{sec:proofofThm1}

We now give a brief outline of the proof of Theorem~\ref{Ramseythm}:
consider any
red/blue colouring of the hyperedges of~$K_m^{(k)}$, where $m=C |\Hy|$
and~$C$ is a large constant depending only on~$k$ and the maximum degree of~$\Hy$.
We apply the hypergraph regularity lemma to the red
subhypergraph~$\G_{red}$ to obtain a reduced
hypergraph~$\mathcal{R}$ which is very dense.
Thus the following fact will show that~$\mathcal{R}$ contains a
copy of~$K_\ell^{(k)}$ with $\ell:=R(K_{k\Delta}^{(k)})$.
% The lemma itself is a weak hypergraph analogue of Tur\'an's theorem and can
% be found in~\cite{hyperturan}.
%
\begin{fact}\label{turanlemma}
For all $\ell,k \in \mathbb{N}$ with $\ell \geq k$, every $k$-uniform
hypergraph~$\mathcal{R}$ on $t \geq \ell$ vertices with
$e(\mathcal{R}) > \left(1-{\ell \choose k}^{-1} \right) \binom{t}{k}$ contains a copy of~$K_\ell^{(k)}$.
\end{fact}
\proof
Let $\mathcal{R}$ be as in the statement of the fact.
Assume for the sake of contradiction that $\mathcal{R}$ is $K_\ell^{(k)}$-free.
Then for each $\ell$-subset $S$ of $V(\mathcal{R})$, we have
$e(\mathcal{R}[S])\leq {\ell \choose k} - 1$. But note that
$$ e(\mathcal{R}) = {t- k \choose \ell - k}^{-1} \sum_S e(\mathcal{R}[S]).$$
Thus $e(\mathcal{R}) \leq {t- k \choose \ell - k}^{-1} {t \choose \ell} \left({\ell \choose k} - 1 \right)$.
Now the observation that ${t- k \choose \ell - k}^{-1} {t \choose \ell} {\ell \choose k} ={t \choose k}$
yields the required contradiction.
\endproof
The copy of~$K_\ell^{(k)}$ in~$\mathcal{R}$ involves~$\ell$ clusters and for
each $k$-tuple of them the red hypergraph~$\G_{red}$ is regular
with respect to almost all of the polyads induced on it. We will
then show that we can find a $(k-1,\ell)$-complex~$\cS$ on these
clusters such that for each $j=2,\dots,k-1$ the restriction of its underlying $j$-uniform
hypergraph~$\cS_j$ to any $(j+1)$-tuple of clusters is a polyad.
Moreover, $\G_{red}$ will be regular with respect to~$\cS_{k-1}$.
By combining $E(\G_{red})\cap \K_k(\cS_{k-1})$ with~$\cS$, we will obtain a regular
$k$-complex~$\cS_{red}$. Similarly we obtain a $k$-complex~$\cS_{blue}$
which also turns out to be regular. We then consider the following red/blue colouring
of~$K_\ell^{(k)}$. We colour a hyperedge red
if~$\G_{red}$ has density at least~$1/2$ with respect to the
corresponding polyad in~$\cS_{k-1}$ and blue otherwise.
By the definition of~$\ell$, we can find a monochromatic
$K_{k\Delta}^{(k)}$. If it is red, then we can apply the embedding
lemma to~$\cS_{red}$ to find a red copy of~$\Hy$. This can be done
since $\Delta(\Hy)\le \Delta$
implies that the chromatic number of~$\Hy$ is at most~$(k-1)\Delta+1 \leq
k\Delta$. If our monochromatic copy of~$K_{k\Delta}^{(k)}$ is blue,
then we can apply the embedding theorem
to~$\cS_{blue}$ and obtain a blue copy of~$\Hy$.

\bigskip

\noindent \textbf{Proof of Theorem \ref{Ramseythm}.}
Given $\Delta$ and $k$, we choose $C$ to be a sufficiently large constant.
We will describe the bounds that~$C$ has to satisfy at the end of the proof.
Let $m:=C|\Hy|$ and consider any red/blue colouring of the hyperedges
of~$K_m^{(k)}$. Let~$\G_{red}$ be the red and~$\G_{blue}$ be the
blue subhypergraph on~$V=V(K_m^{(k)})$. We may assume
without loss of generality that $e(\G_{red})\geq e(\G_{blue})$.%
     \COMMENT{this is not necessary, but saves the reader from having to
check that it is ok to apply the RL to an extremly sparse graph}
We apply the hypergraph regularity lemma to~$\G_{red}$ with
constants~$\eta,\delta_k\ll 1/\Delta,1/k$ as well as  functions~$r$ and~$\delta$
satisfying the hierarchy in~(\ref{reghierarchy}). This gives us clusters
$V_1,\ldots,V_{a_1}$, each of size~$n$ say, together with a $t$-bounded
$(\eta,\delta,\av)$-equitable family of partitions~$\Part=\Part(k-1,\av)$
on~$V$ where $\av=(a_1,\ldots, a_{k-1})$. (Note that by deleting some vertices
of~$\G_{red}$ if necessary we may assume that~$m=|\G_{red}|$ is divisible by~$t!$.)
Since $\eta\ll 1/\Delta,1/k$, condition~(2) in the definition of
an $(\eta,\delta,\av)$-equitable family of partitions implies that
the~$a_1$ which we obtain from the regularity lemma satisfies
$$a_1\ge R(K_{k\Delta}^{(k)})=:\ell.$$
Note that the definition of~$\ell$ involves
a hypergraph Ramsey number whose value is unknown. However, for the argument below
all we need is that this number exists.

Let~$\R$ denote the reduced hypergraph, defined in the previous section.
Proposition~\ref{badktuples} implies that~$\R$ has at least
$(1-\eps)\binom{a_1}{k}$ hyperedges, where
$\eps:=4\sqrt{\delta_k}k!$. Since $\delta_k\ll 1/\Delta,1/k$,
we may assume that $e(\mathcal{R}) \ge (1-\eps)\binom{|\mathcal{R}|}{k}>
\left(1-{\ell \choose k}^{-1} \right) \binom{|\mathcal{R}|}{k}$. Since
$|\mathcal{R}|= a_1 \ge \ell$, this means that we can apply
Fact~\ref{turanlemma} to~$\mathcal{R}$ to obtain a copy
of~$K_\ell^{(k)}$ in~$\mathcal{R}$. Without loss of generality we may
assume that the vertices of this copy are the clusters $V_1,\dots,V_\ell$.

As mentioned above, we now want to find a
$(k-1,\ell)$-complex~$\cS$ on these clusters such that for each $j=2,\dots,k-1$
its underlying $j$-uniform hypergraph~$\cS_j$ is a union
of parts of~$\Part^{(j)}$ and~$\G_{red}$ is regular with respect
to~$\cS_{k-1}$. We construct~$\cS$ inductively starting from the lower
levels. To begin with, for each pair $V_i,V_j$ ($1\leq i < j \leq
\ell$) independently, we choose with probability~$1/a_2$ one of the parts of~$\Part^{(2)}$
induced on $V_i,V_j$. $\cS_2$ will be the union of these parts.
Now suppose that we have chosen~$\cS_{j-1}$ such that
its restriction to any $j$-tuple of clusters forms a
polyad (clearly this is the case for~$\cS_2$).
Now, if~$\hat{P}^{(j-1)}$ is such a polyad, we choose a part of~$\Part^{(j)}$
uniformly at random among the~$a_j$ parts of~$\Part^{(j)}$ that
form~$\K_j(\hat{P}^{(j-1)})$, independently for each $j$-tuple of clusters.
We let~$\cS$ be the $(k-1,\ell)$-complex thus obtained.

We will show that there is some choice of~$\cS$ such that for every $k$-tuple
among the clusters $V_1,\ldots, V_\ell$ the hypergraph~$\G_{red}$ is
$(\delta_k,r)$-regular with respect to the restriction of~$\cS_{k-1}$ to
this $k$-tuple. Note that~$\cS_{k-1}$ restricted to any particular $k$-tuple of
clusters is in fact a polyad selected uniformly at random among
all polyads~$\hat{P}^{(k-1)}$ induced by these~$k$ clusters. Therefore,
since all the $k$-tuples of clusters are fruitful, the definition of a
fruitful $k$-tuple implies that the probability that~$\G_{red}$ has the necessary
regularity is at least
$$1- \sqrt{\delta_k} \binom{\ell}{k} > \frac{1}{2}.$$
The final inequality holds since we may assume that~$\delta_k$
is sufficiently small compared to~$1/\ell$. This shows the existence of a
$(k-1,\ell)$-complex~$\cS$ with the required properties. In what follows,
$P_{\cS}$ will always denote a $(k-1)$-uniform subhypergraph of~$\cS$
induced by~$k$ of the clusters $V_1,\ldots, V_\ell$. So each such~$P_{\cS}$
is a polyad and to each hyperedge of the subhypergraph of~$\R$ induced
by the clusters $V_1,\ldots, V_\ell$ there corresponds such a polyad~$P_{\cS}$.

We now use the densities of~$\G_{red}$ with respect to~$\cS_{k-1}$
to define a red/blue colouring of the~$K_\ell^{(k)}$ which we
found in~$\R$: we colour a hyperedge of this~$K_\ell^{(k)}$ red if
the polyad~$P_{\cS}$ corresponding to this hyperedge satisfies
$d(\G_{red}|P_{\cS}) \ge 1/2$; otherwise we colour it blue. Since
$\ell=R(K_{k\Delta}^{(k)})$, we find a monochromatic copy~$K$
of~$K_{k\Delta}^{(k)}$ in our~$K_\ell^{(k)}$. We now greedily
assign the vertices of~$\Hy$ to the clusters that form the vertex
set of~$K$ in such a way that if~$k$ vertices of~$\Hy$ form a
hyperedge, then they are assigned to~$k$ different clusters. (We
may think of this as a $(k\Delta)$-vertex-colouring of~$\Hy$.) We
now need to show that with this assignment we can apply the
embedding lemma to find a monochromatic copy of~$\Hy$ in either
the subhypergraph of~$\G_{red}$ induced by the~$k\Delta$ clusters
in~$K$ or the subhypergraph of~$\G_{blue}$ induced by these
clusters.

First suppose that~$K$ is red, so we want to apply the embedding
theorem to the $k$-complex formed by~$\G_{red}$ and~$\cS$ (induced
on the~$k\Delta$ clusters in~$K$). However, the embedding theorem
requires all the densities involved to be equal and of the
from~$1/a$ for $a\in\mathbb{N}$, whereas all we know is that for
every polyad~$P_{\cS}$ corresponding to a hyperedge of~$K$, we
have $d(\G_{red}|P_{\cS}) \geq 1/2$. This minor obstacle can be
overcome by choosing a subhypergraph $\G_{red}'\subseteq \G_{red}$
such that~$\G_{red}'$ is $(1/2,3\delta_k,r)$-regular with respect
to each polyad~$P_{\cS}$.
The existence of such a~$\G_{red}'$ follows immediately from the slicing lemma (Lemma~\ref{Slice}).%
    \COMMENT{Comment: It's a bit unfortunate that the proof is
actually not in that paper, but since the result really is quite
easy, it's not worth giving further details here}
%It is easy to see that such a~$\G_{red}'$ exists: fix a $k$-tuple
%$V_{i_1},\ldots, V_{i_k}$ that is a hyperedge of~$K$ and
%the corresponding polyad~$P_{\cS}$. (So $P_{\cS}=\cS_{k-1}[V_{i_1},\ldots, V_{i_k}]$.)
%Consider a random subset of the hyperedges of~$\G_{red}$
%which lie in~$\K_k(P_{\cS})$ such that the
%expected density of this random subset is
%$(1\pm \delta_k)/2$ with respect to~$P_{\cS}$. Then this random subset
%has the desired properties with high probability. A similar result is stated
%explicitly in Proposition~33 of~\cite{roedlschacht}, which
%one can apply directly to obtain the bound of $3\delta_k$ on the regularity
%stated above.
We then add $E(\G_{red}')\cap \K_k(\cS_{k-1})$ to the subcomplex of~$\cS$
induced by the clusters in~$K$ to obtain a regular $(k,k\Delta)$-complex~$\cS_{red}$
and we apply the embedding theorem (Theorem~\ref{embhgph}) there to
find a copy of~$\Hy$ in~$\G_{red}'$, and therefore also in~$\G_{red}$.

On the other hand, if~$K$ is blue, we need to prove
that~$\G_{blue}$ is regular with respect to all chosen polyads~$P_{\cS}$. So
suppose $\textbf{Q}=(Q(1),\ldots,Q(r))$ is an $r$-tuple of subhypergraphs
of one of these polyads~$P_{\cS}$, satisfying $|\K_k(\textbf{Q})|>\delta_k |\K_k
(P_{\cS})|$. Let~$d$ be such that~$\G_{red}$ is $(d,\delta_k,r)$-regular
with respect to~$P_{\cS}$. Then
\begin{align*}|(1-d)-d(\G_{blue}|\textbf{Q} )|
 = |d-(1-d(\G_{blue}|\textbf{Q} ))|
 = |d-d(\G_{red}|\textbf{Q} )|
 < \delta_k.
\end{align*}
Thus~$\G_{blue}$ is $(1-d,\delta_k,r)$-regular with respect to~$P_{\cS}$
(note that $\delta_k \ll  1/2 \le 1-d$). Following the same
argument as in the previous case, we add $E(\G'_{blue}) \cap
\K_k(\cS_{k-1})$ to the subcomplex of~$\cS$
induced by the clusters in~$K$ to derive the regular $(k,k\Delta)$-complex
$\cS_{blue}$ to which we can apply the embedding theorem to obtain a
copy of~$\Hy$ in~$\G_{blue}$.

It remains to check that we can choose~$C$ to be a constant depending only
on~$\Delta$ and~$k$. Note that the constants and functions
$\eta$, $\delta_k$, $r$ and~$\delta$
we defined at the beginning of the proof all depend only on~$\Delta$ and~$k$.
So this is also true for the integers~$m_0$ and~$t$ and the
vector $\av = (a_1,\ldots, a_{k-1})$ which we then obtained from the
regularity lemma. Note that in order to be able to apply the regularity lemma
to~$\G_{red}$ we needed $m \geq m_0$, where $m=C|\Hy|$.
This is certainly true if we set $C \ge m_0$.
The embedding theorem allows us to embed subcomplexes of size at most $cn$,
where $n$ is the cluster size and where $c$ satisfies  $c \ll
1/a_2, ,\dots,1/a_{k-1},d_k, 1/(k\Delta)$
(recall that $d_k=1/2$ and $d_i=1/a_i$ for all $i=2,\dots,k-1$). Thus~$c$ too depends
only on~$\Delta$ and~$k$. In order to apply the embedding theorem we needed that
$n\geq n_0$, where~$n_0$ as defined in the embedding theorem depends
only on~$\Delta$ and~$k$. Since the number of clusters is at most~$t$,
this is satisfied if $m \geq t n_0$, which in turn is certainly true if
$C \ge t n_0$. When we applied the embedding
lemma to $\Hy$, we needed that  $|\Hy|\le cn$. Since
$n= m/a_1=C|\Hy|/a_1\ge C|\Hy|/ t$, it
suffices to choose $C\ge t/c$ for this.
Altogether, this shows that we can define the  constant~$C$ in
Theorem~\ref{Ramseythm} by  $C:=\max \{m_0, t n_0, t/c \}$. \endproof

%%%%%%%%%%%%%%%%%%%%%%%%%%%%%%%%%%%%%%%%%%%%%%%%%%%%%%%%%%%%%%%%%%%%%%%%%%%%%%%%%%%%%%%%%

\section{Deriving Lemmas~\ref{countinglemma} and~\ref{densecounting} from earlier work} \label{sec:CL}

First, we deduce Lemma~\ref{densecounting} from~\cite[Cor.~6.11]{KRS}.
The difference between the two is that the latter result only counts complete hypergraphs
but on the other hand it allows for different densities within each level.
We need a few definitions that make this notion precise.
Let~$\G$ be a $(k,t)$-complex. Recall that~$\G_i$ denotes the underlying $i$-uniform
hypergraph of~$\G$. For each $3\le i< k$, we say that~\emph{$\G_i$ is $(\geq
d_i,\delta_i)$-regular with respect to~$\G_{i-1}$}, if
for every $i$-tuple~$\Lambda_i$ of vertex classes of~$\G$ the induced
hypergraph $\G_i[\Lambda_i]$ is $(d_{\Lambda_i},\delta_i)$-regular with respect to
$\G_{i-1}[\Lambda_i]$, for some $d_{\Lambda_i} \geq d_i$. Similarly we define
when~\emph{$\G_k$ is $(\geq d_k,\delta_k,r)$-regular with respect to~$\G_{k-1}$}
and when \emph{$\G_2$ is $(\geq d_2,\delta_2)$-regular}.
Let $\mathbf{d}:=(d_k,\ldots,d_2)$. We say that a
$(k,t)$-complex $\G$ is \emph{$(\geq \mathbf{d},\delta_k,\delta,r)$-regular} if
\begin{itemize}
\rm \item $\G_k$ is $(\geq d_k,\delta_k,r)$-regular with respect
to $\G_{k-1}$;
\rm \item $\G_i$ is $(\geq d_i,\delta)$-regular with respect to
$\G_{i-1}$ for each $3 \leq  i< k$;
\rm \item $\G_2$ is $(\geq d_2,\delta)$-regular.
\end{itemize}

\begin{lemma}[Dense counting lemma for complete complexes \cite{KRS}] \label{KRSthm}
Let $k,t,n_0$ be positive integers and let
$\eps,d_2,\ldots,d_{k-1},\delta$ be positive constants such that
$$ 1/n_0 \ll \delta \ll \eps \ll d_2,\ldots,d_{k-1},1/t.
$$
Then the following holds for all integers $n\ge n_0$.  Suppose that~$\G$ is a
$(\ge (d_{k-1},\dots,d_2),\delta,\delta,1)$-regular $(k-1,t)$-complex with
vertex classes $V_1,\dots,V_t$, all of size~$n$. Then
$$|K_t^{(k-1)}|_\G = (1 \pm \eps) n^t \prod_{i=2}^{k-1} \prod_{\Lambda_i}d_{\Lambda_i},$$
where the second product is taken over all $i$-tuples~$\Lambda_i$ of vertex classes of~$\G$.
\end{lemma}
We now show how to deduce Lemma~\ref{densecounting} from this.
Full details can be found in~\cite{Olly}.

\medskip

\noindent
{\bf Proof of Lemma~\ref{densecounting}.}
First we prove the lemma for the case when $\ell=t$, i.e.~when each of
the vertex classes $X_1,\ldots, X_t$ of~$\Hy$ consists of exactly one
vertex, say $X_i:=\{h_i\}$. Given such an~$\Hy$ and a complex~$\G$ as in Lemma~\ref{densecounting},
we construct a complex~$\G'$ from~$\G$ as follows: Starting with $i=2$,
for all $i$ with $2\le i \le k-1$
in turn, we successively consider each $i$-tuple
$\Lambda_i=(V_{j_1}, \dots,V_{j_i})$ of vertex classes of~$\G$.
If $h_{j_1},\dots,h_{j_i}$ forms an $i$-edge of~$\Hy$
we let $\G'_i[\Lambda_i]=\G_i[\Lambda_i]$. If $h_{j_1},\dots,h_{j_i}$ does not form an~$i$-edge
we make each copy of~$K_i^{(i-1)}$ in $\G'_{i-1}[\Lambda_i]$ into an $i$-edge
of~$\G_{i}'$. Thus in the latter case
the density of $\G'_i[\Lambda_i]$ with respect to $\G'_{i-1}[\Lambda_i]$ will be~1.
(If $i=2$, this means that we let $\G'_{i}[\Lambda_i]$ be the complete bipartite graph
with vertex classes~$V_{j_1}$ and~$V_{j_2}$.)
Using that~$\Hy$ is a complex, it is easy to see that~$\G'$ is also
$(\ge (d_{k-1},\dots,d_2),\delta,\delta,1)$-regular.
Clearly, there is a bijection between the copies of~$\Hy$ in~$\G$ and the copies
of~$K_t^{(k-1)}$ in~$\G'$. So $|\Hy|_\G = |K_t^{(k-1)}|_{\G'}$.
The result now follows if we apply Lemma~\ref{KRSthm} to~$\G'$.

It now remains to deduce Lemma~\ref{densecounting} for arbitrary $\ell$-partite
complexes~$\Hy$ from the result for the above case. For this,
we use a simple argument that was also used in~\cite{CNKO} to
obtain Lemma~\ref{countinglemma} in the case $k=3$.
We define a complex~$\G^*$ from~$\G$ by making~$|X_i|$ copies
$V_i^{1},\ldots, V_i^{|X_i|}$ of each vertex class~$V_i$
in such a way that for any selection of indices $i_1,\ldots, i_t$ the complex
$\G^*[V_1^{i_1},\ldots, V_t^{i_t}]$ is isomorphic to~$\G$.
Note that~$\G^*$ is $|\Hy|$-partite. Also, we can turn~$\Hy$ into an
$|\Hy|$-partite complex~$\Hy^*$ by viewing each vertex as a single vertex class.
Note that different copies of~$\Hy$ in~$\G$ give rise to
different copies of~$\Hy^*$ in~$\G^*$. Thus $|\Hy|_{\G} \leq |\Hy^*|_{\G^*}$.
Conversely, the only case where a copy of~$\Hy^*$ in~$\G^*$ does not correspond to a
copy of~$\Hy$ in~$\G$ is when there is some~$i$ and indices $j_1 \not = j_2$ such
that the vertices that are used by~$\Hy^*$ in~$V_i^{j_1}$ and~$V_i^{j_2}$
correspond to the same vertex of~$V_i$.
It is easy to see that the number of such copies is comparatively small.
Thus the desired bounds on $|\Hy|_{\G}$ immediately follow from the bounds on
$|\Hy^*|_{\G^*}$ which we obtained in the previous paragraph.
\endproof

We now prove Lemma~\ref{countinglemma}.
Its proof is based on the following version of the counting lemma that
accompanies the hypergraph regularity lemma (Theorem~\ref{reglemma}) from~\cite{roedlschacht}.
%We will state and apply this regularity lemma in Sections~\ref{regularity}
%and~\ref{sec:proofofThm1}.
Theorem~\ref{cntlemRS} gives a lower bound on
the number of complete complexes~$K_t^{(k)}$ in a regular $(k,t)$-complex~$\G$,
under less restrictive assumptions on the regularity constants than those
in Lemma~\ref{KRSthm}.

\begin{theorem}[Counting lemma for complete complexes~\cite{roedlschacht2}]\label{cntlemRS}
Let $k,r,t,n_0$ be positive integers and let $\eps,d_2,\ldots,
d_k,\delta,\delta_k$ be positive constants such that $1/d_i \in\mathbb{N}$
for $i=2,\dots,k-1$ and
$$1/n_0\ll 1/r, \delta \ll\min\{\delta_k,d_2,\ldots,d_{k-1}\}\le\delta_k\ll
\eps,d_k,1/t.$$%
Then the following holds for all integers $n\ge n_0$. Suppose that~$\G$ is a
$(\mathbf{d},\delta_k,\delta,r)$-regular $(k,t)$-complex with
vertex classes $V_1,\dots,V_t$, all of size~$n$, which respects the partition
of~$K_t^{(k)}$. Then
$$|K_t^{(k)}|_\G\geq (1- \eps) n^{t} \prod_{i=2}^k d_i^{{k \choose i}}.$$
\end{theorem}
Lemma~\ref{countinglemma} is more general in the sense that it counts copies
of complexes that may not be complete, and also gives an upper bound on their number.
We will deduce Lemma~\ref{countinglemma} from Theorem~\ref{cntlemRS} in several steps.
The first (and main) step is to deduce a counting lemma
which gives the number of copies of complete complexes, but now
in a $(k,t)$-complex~$\G$ where the density of~$\G_i[\Lambda_i]$ with respect to~$\G_{i-1}[\Lambda_i]$
might be different for different $i$-tuples $\Lambda_i$ of vertex classes of~$\G$.

\begin{lemma}[Counting lemma for complete complexes -- different densities]\label{DifDensitiesCL}
Let $k,r,t,n_0$ be positive integers and let $\eps,d_2,\ldots,
d_k,\delta,\delta_k$ be positive constants such that
$$1/n_0\ll 1/r, \delta \ll\min\{\delta_k,d_2,\ldots,d_{k-1}\}\le\delta_k\ll
\eps, d_k,1/t.$$%
Then the following holds for all integers $n\ge n_0$.  Suppose $\G$ is a
$(\geq\mathbf{d},\delta_k,\delta,r)$-regular $(k,t)$-complex with
vertex classes $V_1,\dots,V_t$, all of size~$n$, such that for all $2< i < k$ and
all $i$-tuples~$\Lambda_i$ of vertex classes of~$\G$ the hypergraph
$\G_i[\Lambda_i]$ is $(d_{\Lambda_i}, \delta)$-regular with respect to
$\G_{i-1}[\Lambda_i]$ where~$d_{\Lambda_i}$ can be written as
$d_{\Lambda_i}=p_{\Lambda_i}/q_{\Lambda_i}$ such that
$p_{\Lambda_i},q_{\Lambda_i}\in \mathbb{N}$ and $1/q_{\Lambda_i}\ge d_i$.
Suppose that the analogue holds for all the~$d_{\Lambda_2}$
and all the~$d_{\Lambda_k}$.
Then
$$|K_t^{(k)}|_\G = (1\pm \eps) n^{t} \prod_{i=2}^k \prod_{\Lambda_i}d_{\Lambda_i},$$
where the second product is taken over all $i$-tuples~$\Lambda_i$ of
vertex classes of~$\G$.
\end{lemma}
\proof
We will first prove the lower bound in this lemma by an inductive argument, in
which we allow for different densities in the top levels but not in the
lower levels, and show that we can always move down another
level, until we allow different densities in all levels.
This leads to the following definition. For any $2 < j \leq k$, we
say that a complex $\G$ is $(\geq d_k,\ldots,\geq
d_j,d_{j-1},\ldots,d_2,\delta_k,\delta,r)$-\emph{regular} if
\begin{itemize}
\rm \item $\G_k$ is $(\geq d_k,\delta_k,r)$-regular with respect
to $\G_{k-1}$;
\rm \item $\G_i$ is $(\geq d_i,\delta)$-regular with respect to
$\G_{i-1}$ for each $j \leq i \leq k-1$;
\rm \item $\G_i$ is $(d_i,\delta)$-regular with respect to
$\G_{i-1}$ for each $3 \leq i \leq j-1$;
\rm \item $\G_2$ is $(d_2,\delta)$-regular. % if $j>2$ and $(\ge d_2,\delta)$-regular if $j=2$.
\end{itemize}
\noindent
Choose new constants $\eta_i,\xi_i, \eps_i$ and integers~$r_i$ satisfying
\begin{align*}
1/n_0 \ll \delta=\xi_2 & \ll \dots \ll \xi_k \ll \min\{\delta_k,d_2,\dots,d_{k-1}\}
\le \delta_k=\eta_2 \ll \dots \ll\eta_{k+1} \\
& \ll \eps_k\ll \dots\ll \eps_2=\eps,d_k,1/t
\end{align*} and
$1/n_0\ll 1/r=1/r_2\ll \dots\ll 1/r_k\ll \min\{\delta_k,d_2,\dots,d_{k-1}\}$.
Then the following claim immediately implies the lemma:

\medskip

\noindent
\textbf{Claim.} \emph{Let $2\leq j \le k$. Suppose that~$\G$ satisfies
the conditions of Lemma~\ref{DifDensitiesCL} but is $(\geq d_k,\ldots,\geq
d_j,d_{j-1},\ldots,d_2,\eta_j,\xi_j,r_j)$-regular instead of
$(\geq\mathbf{d},\delta_k,\delta,r)$-regular if $j>2$, where $1/d_i \in \Nat$ for all $i=2,\ldots, j-1$.
Then
$$|K_t^{(k)}|_\G \ge (1-\eps_j) n^t \left( \prod_{i=2}^{j-1} d_i^{\binom{t}{i}} \right)
\prod_{i=j}^{k} \prod_{\Lambda_i} d_{\Lambda_i} .$$}

\smallskip

\noindent We prove this claim by backward induction on~$j$ as follows:
given a $t$-partite complex $\G$ which is $(\geq d_k,\ldots,\geq
d_j,d_{j-1},\ldots,d_2,\eta_j,\xi_j,r_j)$-regular, we will
partition the hyperedges of~$\G_j$ to obtain several
$(\geq d_k,\ldots,\geq d_{j+1},d'_j,d_{j-1},\ldots,d_2,\eta_{j+1},\xi_{j+1},r_{j+1})$-regular
complexes for some~$d'_j$. We will then apply the lower bound from the induction
hypothesis to each of these complexes. Summing over all of them will give the
lower bound in the claim.

We first consider the case $j=k$. We will apply the slicing
lemma (Lemma~\ref{Slice}) to split the $k$th level~$\G_k$ of the
complex~$\G$ to obtain regular complexes whose densities within the $k$th level are
the same. Set $d'_k :=1/\prod_{\Lambda_k} q_{\Lambda_k}$.
The slicing lemma implies that for all~$\Lambda_k$ there is a
partition~$P(\Lambda_k)$ of the set $E(\G_k[\Lambda_k])$ of $k$-edges induced
on~$\Lambda_k$ such that each part is $(d'_k,\eta_{k+1},r_k)$-regular with
respect to~$\G_{k-1}[\Lambda_k]$. So for each~$\Lambda_k$, $P(\Lambda_k)$ has
$d_{\Lambda_k}/d'_k$ parts. Now for each~$\Lambda_k$, choose one part
from~$P(\Lambda_k)$ and let~$\C_k$ denote the resulting $k$-uniform
$t$-partite hypergraph. Let~$\G^{\C_k}$ denote the $k$-complex obtained
from~$\G$ by replacing~$\G_k$ with~$\C_k$. Then
$$
|K_{t}^{(k)}|_{\G} = \sum_{\C_{k}} |K_t^{(k)}|_{\G^{\C_k}}.
$$
Here the summation is over all possible choices of parts from each of the
$\binom{t}{k}$ partitions~$P(\Lambda_k)$. So the number of summands is
$\prod_{\Lambda_k} d_{\Lambda_k}/d'_k ={d'}_k^{-\binom{t}{k}}\prod_{\Lambda_k} d_{\Lambda_k}$.
Moreover, by Theorem~\ref{cntlemRS} each summand in the above sum can be bounded below:
$$
|K_t^{(k)}|_{\G^{\C_k}} \geq (1-\eps_k) n^t
\left(\prod_{i=2}^{k-1} d_i^{t \choose i} \right){d'}_k^{{t \choose k}}.
$$
Altogether, this implies the claim for $j=k$.

Now suppose that $j<k$ and that the claim holds for $j+1$.
To apply the induction hypothesis, we now need to get equal densities
in the $j$th level. We will achieve this by applying the slicing lemma (Lemma~\ref{Slice}) to this level.
Set $d'_j := 1/\prod_{\Lambda_j} q_{\Lambda_j}$. So $1/d'_j \in \mathbb{N}$.
The slicing lemma implies that for every $j$-tuple~$\Lambda_j$
of vertex classes of~$\G$ there is a partition~$P(\Lambda_j)$
of the set $E(\G_j[\Lambda_j])$ of $j$-edges induced on~$\Lambda_j$ such that each part is
$(d'_j,\xi_{j+1})$-regular with respect to~$\G_{j-1}[\Lambda_j]$.
For each~$\Lambda_j$, the corresponding partition~$P(\Lambda_j)$ will have
$a_{\Lambda_j}:=d_{\Lambda_j}/d'_j$ parts. Now for each~$\Lambda_j$, choose one
part from~$P(\Lambda_j)$ and let~$\C_j$ denote the resulting $j$-uniform
$t$-partite hypergraph. We let~$\G^{\C_j}$ denote the $(k,t)$-complex obtained
from~$\G$ as follows: we replace~$\G_j$ by~$\C_j$ and for each $j<i\le k$
we replace~$\G_i$ with the subhypergraph whose $i$-edges are all those $i$-sets of
vertices that span a~$K_i^{(j)}$ in~$\C_j$.
Thus~$\G^{\C_j}_j$ is $(d'_j,\xi_{j+1})$-regular
with respect to~$\G_{j-1}=\G^{\C_j}_{j-1}$.
However, to apply the induction hypothesis this is not enough. We also need to
prove the following more general assertion.
     \textno For all $i=j,\dots,k$ and any~$\Lambda_i$ the following holds.
If $i=j$ then $\G^{\C_j}_i[\Lambda_i]$ is $(d'_j,\xi_{j+1})$-regular
with respect to~$\G^{\C_j}_{i-1}[\Lambda_i]$.
If $j< i< k$ then $\G^{\C_j}_i[\Lambda_i]$ is
$(d_{\Lambda_i},\xi_{j+1})$-regular with respect to $\G^{\C_j}_{i-1}[\Lambda_i]$.
If $i=k$ then $\G^{\C_j}_i[\Lambda_i]$ is
$(d_{\Lambda_i},\eta_{j+1},r_{j+1})$-regular with respect to
$\G^{\C_j}_{i-1}[\Lambda_i]$ for all but at most %8{t \choose k}
$\sqrt{\eta_{j+1}}\prod_{\Lambda_{j}} a_{\Lambda_{j}}$
hypergraphs~$\C_j$.  &(*)

We will prove~$(*)$ by induction on~$i$.
If $i=j$ then we already know that the assertion is true.
So suppose that $i>j$ and that the claim holds for $i-1$.
We will first consider the case when~$i<k$. The induction hypothesis together
with the dense counting lemma for
complete complexes (Lemma~\ref{KRSthm}) implies that
\begin{equation}\label{eqKCj}
|K_{i}^{(i-1)}|_{\G^{\C_j}_{i-1}[\Lambda_{i}]}
\geq \frac{1}{2} n^i \left( \prod_{\ell=2}^{j-1} d_\ell^{i \choose \ell}\right)
{d'}_j^{\binom{i}{j}} \prod_{s=j+1}^{i-1} \prod_{\Lambda_s\subseteq \Lambda_i} d_{\Lambda_s}.
\end{equation}
Similarly, the assumptions on $\G$ in the claim together with Lemma~\ref{KRSthm} imply
\begin{equation}\label{eqKG}
|K_{i}^{(i-1)}|_{ \G_{i-1}[\Lambda_{i}]} \leq 2n^{i}
\left(\prod_{\ell=2}^{j-1} d_\ell^{{i \choose \ell}}\right)
\prod_{s=j}^{i-1}\prod_{\Lambda_s\subseteq \Lambda_i} d_{\Lambda_s}.
\end{equation}
If we combine these inequalities and use the fact that
$\xi_{j} \ll \xi_{j+1} \ll d_j, 1/k$, we obtain
\begin{equation}\label{eqKCjG}
|K_{i}^{(i-1)}|_{\G^{\C_j}_{i-1}[\Lambda_{i}]}
\ge \sqrt{\xi_{j+1}} |K_{i}^{(i-1)}|_{\G_{i-1}[\Lambda_{i}]}
\ge \frac{ \xi_{j}}{\xi_{j+1}} |K_{i}^{(i-1)}|_{\G_{i-1}[\Lambda_{i}]}.
\end{equation}
In other words, a $\xi_{j+1}$-proportion of copies of~$K_{i}^{(i-1)}$
in $\G^{\C_j}_{i-1}[\Lambda_{i}]$ gives rise to a $\xi_{j}$-proportion of copies in
$\G_{i-1}[\Lambda_{i}]$. Moreover, $\K_i(\G^{\C_j}_{i-1}[\Lambda_i])
\cap E(\G_i^{\C_j}[\Lambda_i] )=
\K_i(\G^{\C_j}_{i-1}[\Lambda_i]) \cap E(\G_i [\Lambda_i] )$ by the definition of~$\G^{\C_j}$
and so $d(\G^{\C_j}_{i}[\Lambda_i]|\G^{\C_j}_{i-1}[\Lambda_i])=
d(\G_{i}[\Lambda_i]|\G^{\C_j}_{i-1}[\Lambda_i])=d_{\Lambda_i}\pm \xi_j$
by~(\ref{eqKCjG}) and the $(d_{\Lambda_{i}},\xi_{j})$-regularity
of $\G_i[\Lambda_i]$ with respect to $\G_{i-1}[\Lambda_i]$.
Thus the $(d_{\Lambda_{i}},\xi_{j+1})$-regularity
of $\G^{\C_j}_i[\Lambda_i]$ with respect to $\G^{\C_j}_{i-1}[\Lambda_i]$
follows from the $(d_{\Lambda_{i}},\xi_{j})$-regularity
of $\G_i[\Lambda_i]$ with respect to $\G_{i-1}[\Lambda_i]$.

But if $i=k$, this might not be true, as~$\eta_{j+1}$ may not be small compared to~$d_{j}$.
However, given a $k$-tuple~$\Lambda_k$ of vertex classes of $\G$, it is true for
most complexes $\G^{\C_j}[\Lambda_k]$.
To see this, given~$\Lambda_k$, let~$\B$ be a
$(k,k)$-complex obtained as follows: For each~$\Lambda_j\subset \Lambda_k$,
choose one part from~$P(\Lambda_j)$ and let~$\B_j$ denote the resulting $j$-uniform
$k$-partite hypergraph. To obtain~$\B$ from~$\G[\Lambda_k]$, we replace~$\G_j[\Lambda_k]$
by~$\B_j$ and for each $j<i\le k$ we replace~$\G_i[\Lambda_k]$ with the
subhypergraph whose $i$-edges are all those $i$-sets of vertices which span a~$K_i^{(j)}$ in~$\B_j$.
Thus there are $\prod_{\Lambda_j\subset \Lambda_k} a_{\Lambda_j}=:A_{\Lambda_k}$
such complexes~$\B$. (Recall that $a_{\Lambda_j}=d_{\Lambda_j}/d_j'$
was the number of parts of the partition $P(\Lambda_j)$.)
Using that~($*$) holds for all $i<k$, similarly as
in (\ref{eqKCj})--(\ref{eqKCjG}) one can show that
\begin{equation}\label{eqKB}
|K_{k}^{(k-1)}|_{\B_{k-1}}\ge
\frac{{d_j' }^{\binom{k}{j}}}{4\prod_{\Lambda_j\subset \Lambda_k} d_{\Lambda_j}}
|K_{k}^{(k-1)}|_{\G_{k-1}[\Lambda_{k}]}
= \frac{|K_{k}^{(k-1)}|_{\G_{k-1}[\Lambda_{k}]}}{4A_{\Lambda_k}}.
\end{equation}
We will now prove the following:
   \textno The underlying $k$-uniform hypergraph~$\B_k$ is not
$(d_{\Lambda_k},\eta_{j+1},r_{j+1})$-regular with respect to~$\B_{k-1}$
for less than $\eta_{j+1} A_{\Lambda_k}$ of the complexes~$\B$. &(**)

If~($**$) is false then we can find $T:=\eta_{j+1} A_{\Lambda_k}/2$ such
complexes $\B^{1},\ldots, \B^{T}$, such that each~$\B^\ell$ has a
$\mathbf{Q}^\ell = (Q_{1}^\ell,\ldots, Q_{r_{j+1}}^\ell)$ satisfying $Q_s^\ell \subseteq \B_{k-1}^\ell$
for all $s=1,\dots,r_{j+1}$ and
$|K_{k}^{(k-1)}|_{\mathbf{Q}^\ell} \geq \eta_{j+1}  |K_{k}^{(k-1)}|_{\B^\ell_{k-1}}$, but either
$d(\B^\ell_k|\mathbf{Q}^\ell ) > d_{\Lambda_k} + \eta_{j+1}$ for each~$\ell$ or
$d(\B^\ell_k|\mathbf{Q}^\ell) < d_{\Lambda_k} - \eta_{j+1}$ for each~$\ell$. We will assume the
latter -- the proof in the former case is similar.
But then let $\mathbf{Q}=(\mathbf{Q}^1,\mathbf{Q}^2,\ldots,\mathbf{Q}^{T})$.
Thus~$\mathbf{Q}$ is a $Tr_{j+1}$-tuple and
$$
|K_{k}^{(k-1)}|_\mathbf{Q}\ge \sum_{\ell=1}^{T} \eta_{j+1}|K_{k}^{(k-1)}|_{\B^\ell_{k-1}}
\stackrel{(\ref{eqKB})}{\geq} \eta_{j} |K_{k}^{(k-1)}|_{\G_{k-1}[\Lambda_{k}]}.
$$
Since we may assume that $Tr_{j+1} \leq r_j$ our assumption on the
regularity of~$\G_k[\Lambda_k]$ with respect to~$\G_{k-1}[\Lambda_k]$
implies that $d(\G_k[\Lambda_k] | \mathbf{Q}) \geq d_{\Lambda_k} - \eta_{j}$. On the
other hand, the definition of~$\B$ implies that
$d(\B^\ell_k | \mathbf{Q}^\ell)=d(\G_k[\Lambda_k]| \mathbf{Q}^\ell)$.
Thus
$d(\G_k[\Lambda_k]|\mathbf{Q})\leq \max_{1\leq\ell\leq T} d(\G_k[\Lambda_k]|\mathbf{Q}^\ell)
= \max_{1 \leq \ell \leq T} d(\B^\ell_k | \mathbf{Q}^\ell) < d_{\Lambda_k} -\eta_{j+1}$.%
\COMMENT{The first inequality holds since
$\frac{a_1+\ldots + a_\ell}{b_1 + \ldots + b_\ell}\leq \max \{ {a_i \over b_i}\}=:m$. Indeed,
$a_i\le m b_i$ and thus $a_1+\dots + a_\ell\le m(b_1 + \ldots + b_\ell)$.}
This is a contradiction, and so~$(**)$ holds.

Note that~$(**)$ implies that for all but at most
${t \choose k}\eta_{j+1}\prod_{\Lambda_{j}} a_{\Lambda_{j}}$
hypergraphs~$\C_j$ the hypergraph~$\G^{\C_j}_k$ is
$(d_{\Lambda_k},\eta_{j+1},r_{j+1})$-regular with respect
to~$\G^{\C_j}_{k-1}$ -- we call these $\C_j$ \emph{nice}.
Since $\eta_{j+1}\ll 1/t$, this completes the proof of~$(*)$.

We are now ready to finish the proof of the induction step of the claim. The
induction
%First suppose that $j<k-1$. Then
\begin{eqnarray*}
|K_t^{(k)}|_\G &\geq & \sum_{\text{nice } \C_j} |K_t^{(k)}|_{\G^{\C_j}}
\geq  (1-\eps_{j+1})\sum_{\text{nice } \C_j} n^t
\left(\prod_{i=2}^{j-1}d_{i}^{t \choose i}\right) {d'}_j^{{t \choose j}}
\prod_{i=j+1}^{k} \prod_{\Lambda_i} d_{\Lambda_i}.
\end{eqnarray*}
The summation is over all possible choices of nice~$\C_j$.
So the number of summands is at least $(1-\sqrt{\eta_{j+1}})\prod_{\Lambda_{j}} a_{\Lambda_{j}}$
and for each~$\Lambda_j$ we have $a_{\Lambda_{j}}d'_j=
d_{\Lambda_j}$.
Since $\eta_{j+1}, \eps_{j+1} \ll \eps_j$, the claim follows and hence the lower bound
in Lemma~\ref{DifDensitiesCL} as well.

It is straightforward to obtain a corresponding upper bound from the lower bound.
The proof is based on an argument that was used
in~\cite{NR03} and later in~\cite{CNKO} to derive a similar upper bound in the case of
3-complexes and thus we only give a sketch of it. A detailed proof can be found in~\cite{Olly}.
Let~$[t]^k$ denote the set of all $k$-subsets of~$[t]=\{1,\dots,t\}$.
Given $S \subseteq [t]^{k}$, we let~$\G^S$ denote the $(k,t)$-complex obtained
from~$\G$ as follows: for each $\{i_1,\ldots, i_k \} \in S$ we replace
the set~$E_k(\G[\Lambda_k])$ of all $k$-edges of~$\G$ induced on $\Lambda_k :=\{V_{i_1},\ldots,
V_{i_k}\}$ by $\K_k(\G_{k-1}[\Lambda_k]) \setminus E_k(\G[\Lambda_k])$.
Thus the density of~$\G_k^S[\Lambda_k]$ with respect to~$\G_{k-1}^S[\Lambda_k]$
is now~$1-d_{\Lambda_k}$. Moreover,
$$ |K_t^{(k-1)}|_{\G_{k-1}} = \sum_{S \subseteq [t]^{k}} |K_t^{(k)}|_{\G^S}.$$
Observe that $|K_t^{(k)}|_{\G} =|K_t^{(k)}|_{\G^{\emptyset}}$ and hence
$$|K_t^{(k)}|_{\G} = |K_t^{(k-1)}|_{\G_{k-1}} -
\sum_{S \subseteq [t]^{k}, S \not = \emptyset} |K_t^{(k)}|_{\G^S}.  $$
Thus, to obtain an upper bound on $|K_t^{(k)}|_{\G}$ all we have to do now is to obtain
an upper bound on $|K_t^{(k-1)}|_{\G_{k-1}}$ and a lower bound on $|K_t^{(k)}|_{\G^S}$,
for every non-empty~$S$.
But the former follows from the dense counting lemma for complete complexes (Lemma~\ref{KRSthm})
and the latter follows from the lower bound in Lemma~\ref{DifDensitiesCL}, which we proved above. (This is why
%in Lemma~\ref{DifDensitiesCL}
we need to allow more general densities than just $1/a$, for $a \in \Nat$.)
\endproof

Lemma~\ref{countinglemma} now follows from Lemma~\ref{DifDensitiesCL} in exactly the same way
as Lemma~\ref{densecounting} followed from Lemma~\ref{KRSthm}.

\section{Proof of the Extension Lemmas~\ref{extensions_count} and~\ref{denseextension}}\label{sec:ext}

We now use Lemma~\ref{countinglemma} to derive
Lemma~\ref{extensions_count} (Lemma~\ref{denseextension} can be
derived in the same way from Lemma~\ref{densecounting}). The proof
idea is similar to that of~\cite[Cor.~14]{roedlschacht2},
\cite[Lemma~6.6]{Gowers2} and~\cite[Lemma~5]{CNKO}. Pick a copy
$H$ of $\Hy$ in $\G$ uniformly at random, and define $X:=\htohp$.
Then $X$ is a random variable. We have
$\mathbb{E}(X)=\frac{1}{|\Hy|_{\G}}\sum_{H\in\G} \htohp =
|\Hy'|_{\G}/|\Hy|_{\G}$. (Here the sum $\sum_{H\in\G}$ is over all
copies of~$\Hy$ in~$\G$.) We pick some constant $\eps$ satisfying
$\delta_k \ll \eps \ll \beta$. By applying the upper bound of the
counting lemma (Lemma~\ref{countinglemma}) to $\Hy$ and the lower
bound to $\Hy'$ we obtain a lower bound for $\mathbb{E}(X)$.
Similarly we obtain an upper bound.
In this way we can easily deduce that%
\begin{equation} \label{expX}
\mathbb{E}(X) = (1 \pm \sqrt{\eps}) \exphtohp .
\end{equation} %
Now consider $\mathbb{E}(X^2)$. We aim to show that its value is
approximately $\exphtohp^2$, and so $X$ has a low variance. Using
Chebyshev's inequality, this will then imply that $X$ is
concentrated around its mean. In other words, only a few copies of
$\Hy$ do not extend to the correct number of copies of~$\Hy'$ in~$\G$.

Observe that $\mathbb{E}(X^2) = \frac{1}{|\Hy|_\G} \sum_{H\in\G}
\htohp^2$. We view $\htohp^2$ as the number of pairs $H_1',H_2'$ of
copies of~$\Hy'$ which extend~$H$. Here the pairs are allowed to
overlap, but we first obtain a rough estimate by insisting that they
intersect precisely in~$H$. So let~$\Hy^*$ be the $(k,\ell)$-complex%
     \COMMENT{Saying that $\Hy^*$ is a $(k,\ell)$-complex makes it at
least clearer what the vertex classes are (which is important since
we are counting only partition-respecting copies
of~$\Hy^*$).}
obtained from two disjoint copies of~$\Hy'$ by identifying them on~$\Hy$.
Thus any copy of~$\Hy^*$ in~$\G$ extending $H$ corresponds to
a pair~$H_1',H_2'$. However, we will later need to take account of
those pairs $H_1',H_2'$ which do not arise from a copy of~$\Hy^*$.
These pairs are exactly those whose intersection is strictly larger
than~$H$.

By applying the counting lemma (Lemma~\ref{countinglemma}) to $\Hy^*$ and to $\Hy$, as before
we obtain\COMMENT{In this case we have no additional error terms
because of non-disjointness in the second inequality since we are
dealing with expectations}%
$$
\frac{1}{|\Hy|_\G} \sum_{H\in\G} \htohs = (1 \pm
\sqrt{\eps})\exphtohs = (1 \pm \sqrt{\eps})\exphtohp^2.
$$%
On the other hand, the number of pairs $H_1',H_2'$ which do not
arise from a copy of~$\Hy^*$ is at most%
     \COMMENT{To see the $(t'-t)^2$ choose $H$ first, then choose a vertex to be
used twice, then choose~$\Hy'$.}
$(t'-t)^2n^{2(t'-t)-1} <
\eps((\prod_{i=2}^k d_i^{e_i(\Hy')-e_i(\Hy)})n^{t'-t})^2 = \eps
\exphtohp^2$. Thus%
\begin{equation}\label{expXsquared}
\frac{1}{|\Hy|_\G} \sum_{H\in\G} \htohp^2 = (1 \pm 2 \sqrt{\eps}) \exphtohp^2.
\end{equation}%
Putting (\ref{expX}) and (\ref{expXsquared}) together, we obtain%
$$var(X) = \mathbb{E}(X^2)-(\mathbb{E}(X))^2 < 5 \sqrt{\eps} \exphtohp ^2.$$%
Now recall Chebyshev's
inequality: $ \mathbb{P}(|X-\mathbb{E}(X)| \geq t) \leq var(X)/t^2.$
We apply this inequality with $t:=\beta \exphtohp$. This implies
that the probability that a randomly chosen copy of~$\Hy$ in~$\G$ does not
satisfy the conclusion of the extension lemma is at most
$var(X)/\beta^2\exphtohp^2 < 5\sqrt{\eps}/\beta^2 < \beta$, and so
at most $\beta |\Hy|_\G$ copies of~$\Hy$ do not satisfy the
conclusion, as required.

\section{Acknowledgement}

We are grateful to the referees for their detailed comments.

\medskip

\noindent
{\footnotesize
Oliver Cooley, Nikolaos Fountoulakis, Daniela K\"{u}hn \& Deryk Osthus,\\
School of Mathematics, University of Birmingham, Edgbaston, Birmingham B15 2TT, UK}

{\footnotesize \parindent=0pt
\it{E-mail addresses}:
\tt{\{cooleyo,nikolaos,kuehn,osthus\}@maths.bham.ac.uk}}

\end{document}